\documentclass[a4paper,12pt]{article}
\usepackage[cp1251]{inputenc} % Для WiN-кодировки в MiKTeX
\usepackage[english]{babel}
\usepackage{amsfonts,amssymb}
\title{Dihedral and reflexive modules with $\infty$-simplicial faces and dihedral and reflexive
homology of involutive $\A$-algebras over unital commutative rings.}
\author{S.V. Lapin}
\date{}
\newcommand{\F}{F_{\infty}}
\newcommand{\D}{D_{\infty}}
\newcommand{\bu}{\bullet}
\newcommand{\E}{E_{\infty}}
\newcommand{\A}{A_{\infty}}

\newcommand{\Hom}{{\rm Hom}}
\newcommand{\p}{\partial}

\begin{document}
\maketitle

\begin{abstract} The concepts of a dihedral
and a reflexive module with $\infty$-simplicial faces are
introduced. For each involutive $A_\infty$-algebra, the dihedral
and the reflexive tensor modules with $\infty$-simplicial faces
are constructed. On the basis of dihedral and reflexive modules
with $\infty$-simplicial faces that defined by an involutive
$A_\infty$-al\-gebra the constructions of the dihedral and the
reflexive homology of involutive $A_\infty$-algebras over any
unital commutative rings are given. The conception of an
involutive homotopy unital $A_\infty$-algebra is introduced. A
long exact sequence that connecting the dihedral and the reflexive
homology of involutive homotopically unital $A_\infty$-algebras
over any unital commutative rings is constructed.
\end{abstract}

Dihedral homology of involutive associative algebras over fields
of characteristic zero was first defined in \cite{Tsygan} as the
homology of the complex of coinvariants of the action of the
dihedral group on the Hochschild complex of these involutive
associative algebras. After that in \cite{KLS} (see also
\cite{KLS1}) the theory of dihedral simplicial modules and in
particular the theory of dihedral modules with simplicial faces
was developed. Further on the basis of the combinatorial technique
of dihedral modules with simplicial faces the dihedral homology
theory of involutive associative algebras over any unital
commutative rings was constructed. In \cite{KLS} also was show
that over fields of characteristic zero the definitions from
\cite{Tsygan} and \cite{KLS} are equivalent. Constructed in
\cite{KLS} the dihedral homology theory proved to be a very useful
tool in the study of the hermitian algebraic $K$-theory of
involutive unital associative algebras, algebraic and homotopy
properties of the hermitian algebraic $K$-theory of topology
spaces  and also in the study of the rational homotopy type of
groups of diffeomorphisms of smooth manifolds (see, g.e., the
survey \cite{Sol}).

On the other hand in \cite{B} the dihedral homology of involutive
$A_\infty$-algebras over fields of characteristic zero was defined
as the homology of the complex of coinvariants of the action of
the dihedral group on the Hochschild complex of these involutive
$A_\infty$-algebras. In this regard gives rise to the important
and very interesting problem of constructing the dihedral homology
theory of involutive $A_\infty$-algebras over any unital
commutative rings, which generalizes developed in \cite{KLS} the
dihedral homology theory of involutive unital associative algebras
over any unital commutative rings. The interest to this problem is
caused, mainly, by the important question about the possible of
constructing  the hermitian algebraic $K$-theory of involutive
homotopy unital $\A$-algebras over any unital commutative rings by
the analogue to how it was done in \cite{Kar} for involutive
unital associative algebras over any unital commutative rings.
Moreover, great interest to constructing of the dihedral homology
theory of involutive $\A$-algebras over any unital commutative
rings is caused also the possible of applying of this theory to
the study of Kontsevich graph-complexes and the cohomology of
moduli spaces by the analogue to how it was done in \cite{PS} for
the cyclic homology of $\A$-algebras.

The present paper is devoted to the solution of the
above-mentioned problem, namely, to constructing on the basis of
the combinatorial technique of dihedral modules with
$\infty$-simplicial faces the dihedral homology theory of
involutive $A_\infty$-algebras over any unital commutative rings.
The paper consists of three paragraphs. In first paragraph, on the
basis of the conceptions of a differential module with
$\infty$-simplicial faces \cite{Lap1}-\cite{Lap7} and a
$\D$-differential module \cite{Lap9}-\cite{Lap17} we introduce the
notions of a dihedral module with $\infty$-simplicial faces and a
reflexive module with $\infty$-simplicial faces. After that, the
notions of the dihedral homology of dihedral modules with
$\infty$-sim\-pli\-cial faces and the reflexive homology of
reflexive modules with $\infty$-simplicial faces are given. In
second paragraph, we construct the dihedral module with
$\infty$-sim\-pli\-cial faces for each involutive
$A_\infty$-al\-geb\-ra over any unital commutative ring (see
Theorem 2.1). Then, we define the dihedral homology of an
involutive $A_\infty$-al\-geb\-ra over any unital commutative ring
as the dihedral homology of the dihedral module with
$\infty$-sim\-pli\-cial faces determined by the given involutive
$A_\infty$-al\-geb\-ra. Since each dihedral module with
$\infty$-sim\-pli\-cial faces can be view as the reflexive module
with $\infty$-sim\-pli\-cial faces, the reflexive homology of
involutive $A_\infty$-al\-geb\-ras over any unital commutative
rings always are defined. Next, we show that over fields of
characteristic zero the definition of the dihedral homology of
involutive $A_\infty$-al\-geb\-ras introduced here is equivalent
to that proposed in \cite{B} (see Corollary 2.1). Also we consider
properties of the reflexive homology of involutive
$A_\infty$-al\-geb\-ras over fields of characteristic zero (see
Corollary 2.2). In the third paragraph, we introduce the
conception of a involutive homotopy unital $\A$-al\-geb\-ra, which
is the involutive analogue of the conception of a homotopy unital
$\A$-al\-geb\-ra \cite{Lyub} (see also \cite{Lap6}). Next, we
construct an exact sequence that connecting the dihedral and the
reflexive homology of involutive homotopy unital $\A$-al\-geb\-ras
over any unital commutative rings (see Theorem 3.1). This exact
sequence generalizes the well-known Krasauskas-Lapin-Solov'ev
exact sequence \cite{KLS} in the dihedral homology theory of
involutive unital associative algebras.

All modules and maps of modules considered in this paper are,
respectively, $K$-modules and $K$-linear maps of modules, where
$K$ is any unital (i.e., with unit) commutative ring.
\vspace{0.5cm}

\centerline{\bf \S\,1. Dihedral and reflexive modules with
$\infty$-simplicial faces}
\vspace{0.5cm}

In what follows, by a bigraded module we mean any bigraded module
$X=\{X_{n,\,m}\}$, $n\geqslant 0$, $m\geqslant 0$, and by a
differential bigraded module, or, briefly, a differential module
$(X,d)$, we mean any bigraded module $X$ endowed with a
differential $d:X_{*,\bu}\to X_{*,\bu-1}$ of bidegree $(0,-1)$.

Recall that a differential module with simplicial faces is defined
as a differential module $(X,d)$ together with a family of module
maps $\p_i:X_{n,\bu}\to X_{n-1,\bu}$, $0\leqslant i\leqslant n$,
which are maps of differential modules and satisfy the simplicial
commutation relations $\p_i\p_j=\p_{j-1}\p_i$, $i<j$. The maps
$\p_i:X_{n,\bu}\to X_{n-1,\bu}$ are called the simplicial face
operators or, more briefly, the simplicial faces of the
differential module $(X,d)$.

Now, we recall the notion of a differential module with
$\infty$-simplicial faces \cite{Lap2} (see also
\cite{Lap3}-\cite{Lap7}), which is a homotopy invariant analogue
of the notion of a differential module with simplicial faces.

Let $\Sigma_k$ be the symmetric group of permutations on a
$k$-element set. Given an arbitrary permutation
$\sigma\in\Sigma_k$ and any $k$-tuple of nonnegative integers
$(i_1,\dots,i_k)$, where $i_1<\dots<i_k$, we consider the
$k$-tuple $(\sigma(i_1),\dots,\sigma(i_k))$, where $\sigma$ acts
on the $k$-tuple $(i_1,\dots,i_k)$ in the standard way, i.e.,
permutes its components. For the $k$-tuple
$(\sigma(i_1),\dots,\sigma(i_k))$, we define a $k$-tuple
$(\widehat{\sigma(i_1)},\dots,\widehat{\sigma(i_k)})$ by the
following formulae $$\widehat{\sigma
(i_s)}=\sigma(i_s)-\alpha(\sigma(i_s)),\quad 1\leqslant s\leqslant
k,$$ where each $\alpha(\sigma(i_s))$ is the number of those
elements of $(\sigma(i_1),\dots,\sigma(i_s),\dots\sigma(i_k))$ on
the right of $\sigma(i_s)$ that are smaller than $\sigma(i_s)$.

A differential module with $\infty$-simplicial faces or, more
briefly, an $\F$-module $(X,d,\widetilde{\p})$ is defined as a
differential module $(X,d)$ together with a family of module maps
$$\widetilde{\p}=\{\p_{(i_1,\dots ,i_k)}:X_{n,\bu}\to
X_{n-k,\bu+k-1}\},\quad 1\leqslant k\leqslant n,$$
$$i_1,\dots,i_k\in\mathbb{Z},\quad 0\leqslant
i_1<\dots<i_k\leqslant n,$$ which satisfy the relations
$$d(\p_{(i_1,\dots,i_k)})=\sum_{\sigma\in\Sigma_k}\sum_{I_{\sigma}}
(-1)^{{\rm sign}(\sigma)+1}
\p_{(\widehat{\sigma(i_1)},\dots,\widehat{\sigma(i_m)})}\,
\p_{(\widehat{\sigma(i_{m+1})},\dots,\widehat{\sigma(i_k)})},\eqno(1.1)$$
where $I_\sigma$ is the set of all partitions of the $k$-tuple
$(\widehat{\sigma(i_1)},\dots,\widehat{\sigma(i_k)})$ into two
tuples $(\widehat{\sigma(i_1)},\dots,\widehat{\sigma(i_m)})$ and
$(\widehat{\sigma(i_{m+1})},\dots,\widehat{\sigma(i_k)})$,
$1\leqslant m\leqslant k-1$, such that the conditions
$\widehat{\sigma(i_1)}<\dots<\widehat{\sigma(i_m)}$ and
$\widehat{\sigma(i_{m+1})}<\dots<\widehat{\sigma(i_k)}$ holds.

The family of maps $\widetilde{\p}=\{\p_{(i_1,\dots ,i_k)}\}$ is
called the $\F$-differential of the $\F$-module
$(X,d,\widetilde{\p})$. The maps $\p_{(i_1,\dots ,i_k)}$ that form
the $\F$-differential of an $\F$-module $(X,d,\widetilde{\p})$ are
called the $\infty$-simplicial faces of this $\F$-module.

It is easy to show that, for $k=1,2,3$, relations $(1.1)$ take,
respectively, the following view $$d(\p_{(i)})=0,\quad i\geqslant
0,\quad d(\p_{(i,j)})=\p_{(j-1)}\p_{(i)}-\p_{(i)}\p_{(j)},\quad
i<j,$$
$$d(\p_{(i_1,i_2,i_3)})=-\p_{(i_1)}\p_{(i_2,i_3)}-\p_{(i_1,i_2)}\p_{(i_3)}-
\p_{(i_3-2)}\p_{(i_1,i_2)}-$$
$$-\,\p_{(i_2-1,i_3-1)}\p_{(i_1)}+\p_{(i_2-1)}\p_{(i_1,i_3)}+\p_{(i_1,i_3-1)}\p_{(i_2)},
\quad i_1<i_2<i_3.$$

Simplest examples of $\F$-modules are differential modules with
simplicial faces. Indeed, given any differential module with
simplicial faces $(X,d,\p_i)$, we can define an $\F$-differential
$\widetilde{\p}=\{\p_{(i_1,\dots ,i_k)}\}:X\to X$ by setting
$\p_{(i)}=\p_i$, $i\geqslant 0$, and $\p_{(i_1,\dots,i_k)}=0$,
$k>1$, thus obtaining the $\F$-module $(X,d,\widetilde{\p})$.

It is worth mentioning that the notion of an $\F$-module specified
above is a part of the general notion of a differential
$\infty$-simplicial module introduced in \cite{Lap3} by using the
homotopy technique of differential Lie modules over curved colored
coalgebras.

Now, we proceed to the notion of a dihedral module with
$\infty$-simplicial faces.

By a dihedral bigraded module $(X,t,r)$ we mean any bigraded
module $X$ together with two families of module maps
$t=\{t_n:X_{n,\bu}\to X_{n,\bu}\}$, $r=\{r_n:X_{n,\bu}\to
X_{n,\bu}\}$, $n\geqslant 0$, satisfying the conditions
$$t_n^{\,n+1}=1_{X_{n,\bu}},\quad r_n^2=1_{X_{n,\bu}},\quad
r_nt_n=t_n^{-1}r_n,\quad n\geqslant 0.$$ In other words, on each
graded module $X_{n,\bu}$, $n\geqslant 0$, the dihedral group of
order $2(n+1)$ with generators $t_n$ and $r_n$ acts on the left.

In what follows, we use the term dihedral differential module for
any quadruple $(X,t,r,d)$, where $(X,t,r)$ is a dihedral bigraded
module, $(X,d)$ is a differential module, and the conditions
$dt_n=t_nd$, $dr_n=r_nd$, $n\geqslant 0$, holds.

Now, recall that a dihedral module with simplicial faces
\cite{KLS} is defined as a dihedral differential module
$(X,t,r,d)$ together with a family of maps $\p_i:X_{n,\bu}\to
X_{n-1,\bu}$, $0\leqslant i\leqslant n$, with respect to which the
triple $(X,d,\p_i)$ is a differential module with simplicial faces
and, moreover, the relations $$\p_it_n=t_{n-1}\p_{i-1},\quad
0<i\leqslant n,\quad \p_0t_n=\p_n,\quad
\p_ir_n=r_{n-1}\p_{n-i},\quad 0\leqslant i\leqslant n,$$ for each
$n\geqslant 0$ are true.

Note that if in the definition of a dihedral module with
simplicial faces we remove the family of automorphisms
$r_n:X_{n,\bu}\to X_{n,\bu}$, $n\geqslant 0$, and the relations
$\p_ir_n=r_{n-1}\p_{n-i}$, $0\leqslant i\leqslant n$, $n\geqslant
0$, then we obtain the definition of a cyclic module with
simplicial faces \cite{Con1}.

{\bf Definition 1.1}. A dihedral module with $\infty$-simplicial
faces or, more briefly, a $D\F$-module, is any five-tuple
$(X,t,r,d,\widetilde{\p})$, where $(X,t,r,d)$ is a dihedral
differential module and $(X,d,\widetilde{\p})$ is a differential
module with $\infty$-simplicial faces related by
$$\p_{(i_1,\dots,i_k)}t_n=\left\{\begin{array}{ll}
t_{n-k}\p_{(i_1-1,\dots,i_k-1)},& i_1>0,\\
(-1)^{k-1}\p_{(i_2-1,\dots,i_k-1,n)},& i_1=0,\\
\end{array}\right.\eqno(1.2)$$
$$\p_{(i_1,\dots,i_k)}r_n=(-1)^{k(k-1)/2}r_{n-k}\p_{(n-i_k,\dots,n-i_1)}.\eqno(1.3)$$

In what follows, we refer to the family of maps
$\widetilde{\p}=\{\p_{(i_1,\dots ,i_k)}:X_{n,\bu}\to
X_{n-k,\bu+k-1}\}$ as the $\F$-dif\-fe\-ren\-tial of the
$D\F$-mo\-du\-le $(X,t,r,d,\widetilde{\p})$. The maps
$\p_{(i_1,\dots ,i_k)}$ that form the $\F$-dif\-fe\-ren\-tial of a
$D\F$-mo\-du\-le $(X,t,r,d,\widetilde{\p})$ are called the
$\infty$-simplicial faces of this $D\F$-module.

Now, we note that if in the definition 1.1 we remove the family of
automorphisms $r_n:X_{n,\bu}\to X_{n,\bu}$, $n\geqslant 0$, and
the relations $(1.3)$ then we obtain the definition of a cyclic
module with $\infty$-simplicial faces \cite{Lap8} or, more
briefly, $C\F$-module. Therefore each $D\F$-module
$(X,t,r,d,\widetilde{\p})$ always defines the $C\F$-module
$(X,t,d,\widetilde{\p})$.

Simple examples of $D\F$-modules are dihedral modules with
simplicial faces. Indeed, given any dihedral module with
simplicial faces $(X,t,r,d,\p_i)$, we can define an
$\F$-dif\-ferential $\widetilde{\p}=\{\p_{(i_1,\dots ,i_k)}\}$ by
setting $\p_{(i)}=\p_i$, $i\geqslant 0$, and
$\p_{(i_1,\dots,i_k)}=0$, $k>1$, thus obtaining the $D\F$-module
$(X,t,r,d,\widetilde{\p})$.

Now, we make preparations to introduce the notion of the dihedral
homology of a dihedral module with $\infty$-simplicial faces.

First, recall that a $\D$-differential module \cite{Lap9} (sees
also \cite{Lap10}-\cite{Lap17}) or, more briefly, a $\D$-module
$(X,d^{\,i})$ is defined as a module $X$ together with a family of
module maps $\{d^{\,i}:X\to X~|~i\in\mathbb{Z},~i\geqslant 0\}$
satisfying the relations
$$\sum\limits_{i+j=k}d^{\,i}d^{\,j}=0,\quad k\geqslant
0.\eqno(1.4)$$ It is worth noting that a $\D$-module $(X,d^{\,i})$
can be equipped with any $\mathbb{Z}^{\times n}$-grading, i.e.,
$X=\{X_{k_1,\dots,k_n}\}$, $(k_1,\dots,k_n)\in\mathbb{Z}^{\times
n}$, $n\geqslant 1$, and maps $d^{\,i}:X\to X$, $i\geqslant 0$,
can have any $n$-degree
$(l_1(i),\dots,l_n(i))\in\mathbb{Z}^{\times n}$, i.e.,
$d^{\,i}:X_{k_1,\dots,k_n}\to X_{k_1+l_1(i),\dots,k_n+l_n(i)}$.
For $k=0$, the relations $(1.4)$ have the form $d^{\,0}d^{\,0}=0$,
and hence $(X,d^{\,0})$ is a differential module. In \cite{Lap9}
was established the homotopy invariance of the structure of a
$\D$-differential module over any unital commutative ring under
homotopy equivalences of differential modules. Later, it was shown
in \cite{LV} that the homotopy invariance of the structure of the
module over fields of characteristic zero can be established by
using the Koszul duality theory.

A $\D$-module $(X,d^{\,i})$ is said to be stable if, for each
$x\in X$, there exists a number $k=k(x)\geqslant 0$ such that
$d^{\,i}(x)=0$, $i>k$. Any stable $\D$-module $(X,d^{\,i})$
determines the differential $\overline{d\,}:X\to X$ defined by
$\overline{d\,}=(d^{\,0}+d^{\,1}+\dots+d^{\,i}+\dots)$. The map
$\overline{d\,}:X\to X$ is indeed a differential, because
relations $(1.4)$ imply the equality
$\overline{d\,}\,\overline{d\,}=0$. It is easy to see that if the
stable $\D$-module $(X,d^{\,i})$ is equipped with a
$\mathbb{Z}^{\times n}$-grading $X=\{X_{k_1,\dots,k_n}\}$,
$k_1\geqslant 0,\dots,k_n\geqslant 0$, and maps $d^{\,i}:X\to X$,
$i\geqslant 0$, have $n$-degree $(l_1(i),\dots,l_n(i))$ satisfying
the condition $l_1(i)+\dots+l_n(i)=-1$, then there is the chain
complex $(\overline{X},\overline{d\,})$ defined by the following
formulae:
$$\overline{X}_m=\bigoplus_{k_1+\dots+k_n=m}X_{k_1\dots,k_n},\quad
\overline{d\,}=\sum_{i=0}^\infty d^{\,i}:\overline{X}_m\to
\overline{X}_{m-1},\quad m\geqslant 0.$$

It was shown in \cite{Lap2} that any $\F$-module
$(X,d,\widetilde{\p})$ determines the sequence of stable
$\D$-modules $(X,d_q^{\,i})$, $q\geqslant 0$, equipped with the
bigrading $X=\{X_{n,m}\}$, $n\geqslant 0$, $m\geqslant 0$, and
defined by the following formulae: $$d_q^{\,0}=d,~~
d_q^{\,k}=\sum\limits_{0\leqslant i_1<\dots<i_k\leqslant n-q
}(-1)^{i_1+\dots+i_k}\p_{(i_1,\dots,i_k)}:X_{n,\bu}\to
X_{n-k,\bu+k-1},~~k\geqslant 1.\eqno(1.5)$$

Let us recall \cite{Lap8} the construction of the chain bicomplex
$(C(\overline{X}),\delta_1,\delta_2)$ that is defined by the
cyclic module with $\infty$-simplicial faces
$(X,t,d,\widetilde{\p})$. Given any $C\F$-mo\-dule
$(X,t,d,\widetilde{\p})$, consider the two $\D$-mo\-dules
$(X,d_0^{\,i})$ and $(X,d_1^{\,i})$ defined by $(1.5)$ for
$q=0,1$, and the two families of maps $$T_n=(-1)^nt_n:X_{n,\bu}\to
X_{n,\bu},\quad n\geqslant 0,$$
$$N_n=1+T_n+T_n^2+\dots+T_n^n:X_{n,\bu}\to X_{n,\bu},\quad
n\geqslant 0.$$ Obviously, the condition $t_n^{\,n+1}=1$,
$n\geqslant 0$, implies the relations $$(1-T_n)N_n=0,\quad
N_n(1-T_n)=0,\quad n\geqslant 0.\eqno(1.6)$$ Moreover, in
\cite{Lap8} it was shown that the families of module maps
$\{T_n:X_{n,\bu}\to X_{n,\bu}\}$, $\{N_n:X_{n,\bu}\to
X_{n,\bu}\}$, $\{d_0^{\,i}:X_{*,\bu}\to X_{*-i,\bu+i-1}\}$ and
$\{d_1^{\,i}:X_{*,\bu}\to X_{*-i,\bu+i-1}\}$ are related by
$$d_0^{\,i}(1-T_n)=(1-T_{n-i})d_1^{\,i},\quad
d_1^{\,i}N_n=N_{n-i}d_0^{\,i},\quad i\geqslant 0,\quad n\geqslant
0.\eqno(1.7)$$ For example, for $i=2$ and $n=3$, we have the
following equalities:
$$d_0^{\,2}(1-T_3)=d_0^{\,2}(1+t_3)=(-\p_{(0,1)}+\p_{(0,2)}-
\p_{(0,3)}-\p_{(1,2)}+\p_{(1,3)}-\p_{(2,3)})(1+t_3)=$$
$$=-\p_{(0,1)}+\p_{(0,2)}-
\p_{(0,3)}-\p_{(1,2)}+\p_{(1,3)}-\p_{(2,3)}-\p_{(0,1)}t_3+\p_{(0,2)}t_3-
\p_{(0,3)}t_3-$$
$$-\p_{(1,2)}t_3+\p_{(1,3)}t_3-\p_{(2,3)}t_3=-\p_{(0,1)}+\p_{(0,2)}-
\p_{(0,3)}-\p_{(1,2)}+\p_{(1,3)}-$$
$$-\p_{(2,3)}-(-1)^{2-1}\p_{(0,3)}+(-1)^{2-1}\p_{(1,3)}-(-1)^{2-1}\p_{(2,3)}-
t_1\p_{(0,1)}+t_1\p_{(0,2)}-t_1\p_{(1,2)}=$$
$$=(1+t_1)(-\p_{(0,1)}+\p_{(0,2)}-\p_{(1,2)})=(1+t_1)d_1^{\,2}=(1-T_1)d_1^{\,2}.$$

$$d_1^{\,2}N_3=d_1^{\,2}(1+T_3+T_3^2+T_3^3)=$$
$$=d_1^{\,2}(1-t_3+t_3^2-t_3^3)=(-\p_{(0,1)}+\p_{(0,2)}-\p_{(1,2)})(1-t_3+t_3^2-t_3^3)=$$
$$=(-\p_{(0,1)}+\p_{(0,2)}-\p_{(1,2)})-
(-\p_{(0,1)}+\p_{(0,2)}-\p_{(1,2)})t_3+(-\p_{(0,1)}+\p_{(0,2)}-\p_{(1,2)})t_3^2-$$
$$-(-\p_{(0,1)}+\p_{(0,2)}-\p_{(1,2)})t_3^3=(-\p_{(0,1)}+\p_{(0,2)}-\p_{(1,2)})-(\p_{(0,3)}-\p_{(1,3)}-t_1\p_{(0,1)})+$$
$$+(-\p_{(2,3)}-t_1\p_{(0,2)}+t_1\p_{(0,3)})-(-t_1\p_{(1,2)}+t_1\p_{(1,3)}-t_1\p_{(2,3)})=$$
$$=(1-t_1)(-\p_{(0,1)}+\p_{(0,2)}-
\p_{(0,3)}-\p_{(1,2)}+\p_{(1,3)}-\p_{(2,3)})=(1-t_1)d_0^{\,2}=$$
$$=(1+T_1)d_0^{\,2}=N_1d_0^{\,2}.$$

Now, we consider the chain complexes $(\overline{X},b)$ and
$(\overline{X},b^{'})$ corresponding to the $\D$-mo\-dules
$(X,d_0^{\,i})$ and $(X,d_1^{\,i})$ specified above; here
$$\overline{X}_n=\bigoplus_{k=0}^n X_{k,n-k},\quad
b=\overline{d\,}_0=\sum_{i=0}^n d_0^{\,i}:\overline{X}_n\to
\overline{X}_{n-1},$$ $$b^{'}=\overline{d\,}_1=\sum_{i=0}^n
d_1^{\,i}:\overline{X}_n\to \overline{X}_{n-1},\quad n\geqslant
0.$$ Consider also the two families of maps
$$\overline{T}_n=\sum_{k=0}^nT_{k}:\overline{X}_n\to\overline{X}_n,\quad
\overline{N}_n=\sum_{k=0}^nN_{k}:\overline{X}_n\to\overline{X}_n,\quad
n\geqslant 0.$$ It is seen from $(1.6)$ and $(1.7)$ that
$$(1-\overline{T}_n)\overline{N}_n=0,\quad
\overline{N}_n(1-\overline{T}_n)=0,\quad n\geqslant 0,$$
$$b(1-\overline{T}_n)=(1-\overline{T}_{n-1})b^{'},\quad
b^{'}\overline{N}_n=\overline{N}_{n-1}b,\quad n\geqslant 0.$$ It
follows from these relations that any $C\F$-module
$(X,t,d,\widetilde{\p})$ determines the chain bicomplex
\vspace{1cm}

\begin{center}
\parbox {4cm}
{\setlength{\unitlength}{1cm}
\begin{picture}(2,3.5)

\put(-3.1,3.5){\makebox(0,0){\vspace{0.3cm}$\vdots$}}
\put(-0.7,3.5){\makebox(0,0){\vspace{0.3cm}$\vdots$}}
\put(1.7,3.5){\makebox(0,0){\vspace{0.3cm}$\vdots$}}
\put(4.1,3.5){\makebox(0,0){\vspace{0.3cm}$\vdots$}}

\put(-3.1,3.2){\vector(0,-1){0.7}}
\put(-0.7,3.2){\vector(0,-1){0.7}}
\put(1.7,3.2){\vector(0,-1){0.7}}
\put(4.1,3.2){\vector(0,-1){0.7}}

\put(-3,2){\makebox(0,0){$\overline{X}_{n+1}$}}
\put(-0.6,2){\makebox(0,0){$\overline{X}_{n+1}$}}
\put(1.8,2){\makebox(0,0){$\overline{X}_{n+1}$}}
\put(4.2,2){\makebox(0,0){$\overline{X}_{n+1}$}}
\put(6.3,2){\makebox(0,0){$\dots$}}

\put(-1.3,2){\vector(-1,0){1.1}} \put(1.1,2){\vector(-1,0){1.1}}
\put(3.5,2){\vector(-1,0){1.1}} \put(5.9,2){\vector(-1,0){1.1}}

\put(-3.3,2.8){\makebox(0,0){$^{b}$}}
\put(-3.3,1.1){\makebox(0,0){$^{b}$}}
\put(-3.3,-0.6){\makebox(0,0){$^{b}$}}
\put(-3.3,-2.3){\makebox(0,0){$^{b}$}}

\put(-1.0,2.8){\makebox(0,0){$^{-b^{'}}$}}
\put(-1.0,1.1){\makebox(0,0){$^{-b^{'}}$}}
\put(-1.0,-0.6){\makebox(0,0){$^{-b^{'}}$}}
\put(-1.0,-2.3){\makebox(0,0){$^{-b^{'}}$}}

\put(1.5,2.8){\makebox(0,0){$^{b}$}}
\put(1.5,1.1){\makebox(0,0){$^{b}$}}
\put(1.5,-0.6){\makebox(0,0){$^{b}$}}
\put(1.5,-2.3){\makebox(0,0){$^{b}$}}

\put(3.8,2.8){\makebox(0,0){$^{-b^{'}}$}}
\put(3.8,1.1){\makebox(0,0){$^{-b^{'}}$}}
\put(3.8,-0.6){\makebox(0,0){$^{-b^{'}}$}}
\put(3.8,-2.3){\makebox(0,0){$^{-b^{'}}$}}

\put(-3.1,1.7){\vector(0,-1){1.0}}
\put(-0.7,1.7){\vector(0,-1){1.0}}
\put(1.7,1.7){\vector(0,-1){1.0}}
\put(4.1,1.7){\vector(0,-1){1.0}}

\put(-3,0.2){\makebox(0,0){$\overline{X}_n$}}
\put(-0.6,0.2){\makebox(0,0){$\overline{X}_n$}}
\put(1.8,0.2){\makebox(0,0){$\overline{X}_n$}}
\put(4.2,0.2){\makebox(0,0){$\overline{X}_n$}}
\put(6.3,0.2){\makebox(0,0){$\dots$}}

\put(-1.3,0.2){\vector(-1,0){1.1}}
\put(1.1,0.2){\vector(-1,0){1.1}}
\put(3.5,0.2){\vector(-1,0){1.1}}
\put(5.9,0.2){\vector(-1,0){1.1}}

\put(-3.1,-0.1){\vector(0,-1){1.0}}
\put(-0.7,-0.1){\vector(0,-1){1.0}}
\put(1.7,-0.1){\vector(0,-1){1.0}}
\put(4.1,-0.1){\vector(0,-1){1.0}}

\put(-3,-1.55){\makebox(0,0){$\overline{X}_{n-1}$}}
\put(-0.6,-1.55){\makebox(0,0){$\overline{X}_{n-1}$}}
\put(1.8,-1.55){\makebox(0,0){$\overline{X}_{n-1}$}}
\put(4.2,-1.55){\makebox(0,0){$\overline{X}_{n-1}$}}
\put(6.3,-1.55){\makebox(0,0){$\dots$}}

\put(-1.7,-1.3){\makebox(0,0){$^{1-\overline{T}_{n-1}}$}}
\put(-1.7,0.4){\makebox(0,0){$^{1-\overline{T}_n}$}}
\put(-1.7,2.2){\makebox(0,0){$^{1-\overline{T}_{n+1}}$}}

\put(0.6,-1.3){\makebox(0,0){$^{\overline{N}_{n-1}}$}}
\put(0.6,0.4){\makebox(0,0){$^{\overline{N}_n}$}}
\put(0.6,2.2){\makebox(0,0){$^{\overline{N}_{n+1}}$}}

\put(3.1,-1.3){\makebox(0,0){$^{1-\overline{T}_{n-1}}$}}
\put(3.1,0.4){\makebox(0,0){$^{1-\overline{T}_n}$}}
\put(3.1,2.2){\makebox(0,0){$^{1-\overline{T}_{n+1}}$}}

\put(5.4,-1.3){\makebox(0,0){$^{\overline{N}_{n-1}}$}}
\put(5.4,0.4){\makebox(0,0){$^{\overline{N}_n}$}}
\put(5.4,2.2){\makebox(0,0){$^{\overline{N}_{n+1}}$}}

\put(-3.1,-2.8){\makebox(0,0){$\vdots$}}
\put(-0.7,-2.8){\makebox(0,0){$\vdots$}}
\put(1.7,-2.8){\makebox(0,0){$\vdots$}}
\put(4.1,-2.8){\makebox(0,0){$\vdots$}}

\put(-1.3,-1.55){\vector(-1,0){1.1}}
\put(1.1,-1.55){\vector(-1,0){1.1}}
\put(3.5,-1.55){\vector(-1,0){1.1}}
\put(5.9,-1.55){\vector(-1,0){1.1}}

\put(-3.1,-1.9){\vector(0,-1){0.7}}
\put(-0.7,-1.9){\vector(0,-1){0.7}}
\put(1.7,-1.9){\vector(0,-1){0.7}}
\put(4.1,-1.9){\vector(0,-1){0.7}}
\end{picture}}
\end{center}
\vspace{3.5cm}

\noindent We denote this chain bicomplex by
$(C(\overline{X}),\delta_1,\delta_2)$, where
$C(\overline{X})_{n,m}=\overline{X}_n$, $n\geqslant 0$,
$m\geqslant 0$, $\delta_1:C(\overline{X})_{n,m}\to
C(\overline{X})_{n-1,m}$, $\delta_2:C(\overline{X})_{n,m}\to
C(\overline{X})_{n,m-1}$, and $$\delta_1=\left\{\begin{array}{ll}
b,&m\equiv 0\,{\rm mod}(2),\\ -b^{'},&m\equiv 1\,{\rm mod}(2),\\
\end{array}\right.\quad\delta_2=\left\{\begin{array}{ll}
1-\overline{T}_n,&m\equiv 1\,{\rm mod}(2),\\
\overline{N}_n,&m\equiv 0\,{\rm mod}(2).\\
\end{array}\right.$$
For the chain complex associated with the chain bicomplex
$(C(\overline{X}),\delta_1,\delta_2)$, we use the notation $({\rm
Tot}(C(\overline{X})),\delta)$, where $\delta=\delta_1+\delta_2$.

Recall \cite{Lap8} that the cyclic homology $HC(X)$ of a
$C\F$-module $(X,t,d,\widetilde{\p})$ is defined as the homology
of the chain complex $({\rm Tot}(C(\overline{X})),\delta)$
associated with the chain bicomplex
$(C(\overline{X}),\delta_1,\delta_2)$.

Given any $D\F$-module $(X,t,r,d,\widetilde{\p})$, consider the
specified above $\D$-modules $(X,d_0^{\,i})$ and $(X,d_1^{\,i})$,
and consider the family of maps
$$R_n=(-1)^{n(n+1)/2}r_n:X_{n,\bu}\to X_{n,\bu},\quad n\geqslant
0.$$ It is easily verified that the relations
$t^{\,n+1}_n=r_n^2=1$, $r_nt_n=t_n^{-1}r_n$, $n\geqslant 0$,
implies the equalities $$(1-T_n)(R_nT_n)=-R_n(1-T_n),\quad
N_nR_n=(R_nT_n)N_n,\quad n\geqslant 0.\eqno(1.8)$$ For any
collection $0\leqslant i_1<\dots<i_k\leqslant n$, the relations
$(1.3)$ implies the equality
$$(-1)^{i_1+\dots+i_k}\p_{(i_1,\dots,i_k)}R_n=
(-1)^{(n-i_k)+\dots+(n-i_1)}R_{n-k}\p_{(n-i_k,\dots,n-i_1)}.$$
Moreover, for any collection $0\leqslant i_1<\dots<i_k\leqslant
n-1$, by using the relations $(1.2)$ and $(1.3)$ we obtain the
equality
$$(-1)^{i_1+\dots+i_k}\p_{(i_1,\dots,i_k)}(R_nT_n)=(-1)^{(n-i_k-1)+\dots+(n-i_1-1)}(R_{n-k}T_{n-k})
\p_{(n-i_k-1,\dots,n-i_1-1)}.$$ For specified above families
$\{d_0^{\,i}:X_{*,\bu}\to X_{*-i,\bu+i-1}\}$ and
$\{d_1^{\,i}:X_{*,\bu}\to X_{*-i,\bu+i-1}\}$, the last two
equalities implies the relations
$$d_0^{\,i}R_n=R_{n-i}d_0^{\,i},\quad
d_1^{\,i}(R_nT_n)=(R_{n-i}T_{n-i})d_1^{\,i},\quad i\geqslant
0,\quad n\geqslant 0.\eqno(1.9)$$ Consider the chain complexes
$(\overline{X},b)$ and $(\overline{X},b^{'})$ that corresponds to
the $\D$-modules $(X,d_0^{\,i})$ and $(X,d_1^{\,i})$. Moreover,
consider the family of maps
$$\overline{R}_n=\sum_{m=0}^nR_m:\overline{X}_n\to\overline{X}_n,\quad
n\geqslant 0.$$ By using the formulae $(1.8)$ and $(1.9)$ we
obtain the following equalities:
$$\left.\begin{array}{c}(1-\overline{T}_n)(\overline{R}_n\overline{T}_n)=-\overline{R}_n(1-\overline{T}_n),
\quad\overline{N}_n\overline{R}_n=(\overline{R}_n\overline{T}_n)\overline{N}_n,\quad
n\geqslant 0,\\
\\
b\overline{R}_n=\overline{R}_{n-1}b,\quad
b^{'}(\overline{R}_n\overline{T}_n)=(\overline{R}_{n-1}\overline{T}_{n-1})b^{'},\quad
n\geqslant 0.\\
\end{array}\right\}\eqno(1.10)$$
Since the $D\F$-module $(X,t,r,d,\widetilde{\p})$ always defines
the $C\F$-module $(X,t,d,\widetilde{\p})$, for the $D\F$-module
$(X,d,t,r,\widetilde{\p})$, the chain bicomplex
$(C(\overline{X}),\delta_1,\delta_2)$ always is defined. The
equalities $(1.10)$ say us that there is a left action of the
group $\mathbb{Z}_2=\{1,\vartheta\,|\,\vartheta^2=1\}$ on the
chain bicomplex $(C(\overline{X}),\delta_1,\delta_2)$ of any
$D\F$-module $(X,t,r,d,\widetilde{\p})$. This left action is
defined by means of the automorphism
$\vartheta:C(\overline{X})_{*,\bu}\to C(\overline{X})_{*,\bu}$ of
the order two, which at any element $x\in C(\overline{X})_{n,m}$
is given by the following rule:
$$\vartheta(x)=\left\{\begin{array}{ll}
(-1)^k\overline{R}_n(x),&m=2k,\\
(-1)^{k+1}\overline{R}_n\overline{T}_n(x),&m=2k+1.\\
\end{array}\right.$$

{\bf Definition 1.2}. We define the dihedral homology $HD(X)$ of a
dihedral module with $\infty$-simplicial faces
$(X,t,r,d,\widetilde{\p})$ as the hyperhomology
$H(\mathbb{Z}_2;(C(\overline{X}),\delta_1,\delta_2))$ of the group
$\mathbb{Z}_2$ with coefficients in
$(C(\overline{X}),\delta_1,\delta_2)$ relative to the specified
above action of the group $\mathbb{Z}_2$ on the chain bicomplex
$(C(\overline{X}),\delta_1,\delta_2)$.

The hyperhomology
$H(\mathbb{Z}_2;(C(\overline{X}),\delta_1,\delta_2))$ is defined
as the homology of the chain complex that associated with the
triple chain complex $({\cal
P}(\mathbb{Z}_2)\otimes_{K[\mathbb{Z}_2]}C(\overline{X}),\delta_1,\delta_2,\delta_3)$,
where $K[\mathbb{Z}_2]$ is a groups algebra of the group
$\mathbb{Z}_2=\{1,\vartheta\,|\,\vartheta^2=1\}$, the chain
complex $({\cal P}(\mathbb{Z}_2),d)$ is any projective resolvente
of the trivial $K[\mathbb{Z}_2]$-module $K$, and the differential
$\delta_3$ is defined by $$\delta_3=(-1)^{n+m}d\otimes 1:{\cal
P}(\mathbb{Z}_2)_l\otimes_{K[\mathbb{Z}_2]}C(\overline{X})_{n,m}\to
{\cal
P}(\mathbb{Z}_2)_{l-1}\otimes_{K[\mathbb{Z}_2]}C(\overline{X})_{n,m}.$$
If we take as the projective resolvente $({\cal
P}(\mathbb{Z}_2),d)$ the standard free resolvente $$({\cal
S}(\mathbb{Z}_2),d):\,K[\mathbb{Z}_2]\buildrel{~1-\vartheta}\over\longleftarrow
K[\mathbb{Z}_2]\buildrel{~1+\vartheta}\over\longleftarrow
K[\mathbb{Z}_2]\buildrel{~1-\vartheta}\over\longleftarrow
K[\mathbb{Z}_2]\buildrel{~1+\vartheta}\over\longleftarrow\cdots,$$
then we obtain that the dihedral homology
$HD(X)=H(\mathbb{Z}_2;(C(\overline{X}),\delta_1,\delta_2))$ of any
$D\F$-module $(X,t,r,d,\widetilde{\p})$ is the homology of the
chain complex associated with the triple chain complex
$(D(\overline{X}),\delta_1,\delta_2,\delta_3)$, where
$D(\overline{X})_{n,m,l}=C(\overline{X})_{n,m}=\overline{X}_n$,
$n\geqslant 0$, $m\geqslant 0$, $l\geqslant 0$, the differentials
$$\delta_1:D(\overline{X})_{n,m,l}=C(\overline{X})_{n,m}\to
C(\overline{X})_{n,m-1}=D(\overline{X})_{n,m-1,l},$$
$$\delta_2:D(\overline{X})_{n,m,l}=C(\overline{X})_{n,m}\to
C(\overline{X})_{n-1,m}=D(\overline{X})_{n-1,m,l}$$ were defined
above, and the differential $\delta_3:D(\overline{X})_{n,m,l}\to
D(\overline{X})_{n,m,l-1}$ is given by
$$\delta_3=\left\{\begin{array}{ll}
(-1)^n(1+(-1)^l\overline{R}_n),&m\equiv 0\,{\rm mod}(4),\\
(-1)^{n+1}(1+(-1)^{l+1}\overline{R}_n\overline{T}_n),&m\equiv
1\,{\rm mod}(4),\\ (-1)^n(1+(-1)^{l+1}\overline{R}_n),&m\equiv
2\,{\rm mod}(4),\\
(-1)^{n+1}(1+(-1)^l\overline{R}_n\overline{T}_n),&m\equiv 3\,{\rm
mod}(4).\\
\end{array}\right.$$
The chain complex associated with the triple chain complex
$(D(\overline{X}),\delta_1,\delta_2,\delta_3)$ we denote by $({\rm
Tot}(D(\overline{X})),\widehat{\delta})$, where
$\widehat{\delta}=\delta_1+\delta_2+\delta_3$.

Note that if a $D\F$-module $(X,t,r,d,\widetilde{\p})$ is a
dihedral module with simplicial faces $(X,t,r,d,\p_i)$, then the
triple chain complex
$(D(\overline{X}),\delta_1,\delta_2,\delta_3)$ coincides with
well-known \cite{KLS} the triple chain complex, which computes the
dihedral homology of the dihedral module with simplicial faces
$(X,t,r,d,\p_i)$.

Now, let us describe the convenient method of a computation of the
dihedral homology of dihedral modules with $\infty$-simplicial
faces over fields of characteristic zero.

Given an arbitrary $D\F$-module $(X,t,r,d,\widetilde{\p})$,
consider the corresponding chain complexes $(\overline{X},b)$ and
$(\overline{X},-b^{'})$. It easily follows from the mentioned
above relations $b(1-\overline{T}_n)=(1-\overline{T}_{n-1})b^{'}$,
$b\overline{R}_n=\overline{R}_{n-1}b$, $n\geqslant 0$, that the
chain complex $(M(\overline{X}),b)$, where
$$M(\overline{X})_n=\overline{X}_n\left/\left({\rm
Im}(1-\overline{T}_n)+{\rm Im}(1-\overline{R}_n)\right),\quad
n\geqslant 0,\right.$$ is well defined.

The following assertion describes the convenient method of a
computation of the dihedral homology of dihedral modules with
$\infty$-simplicial faces over fields of characteristic zero.

{\bf Theorem 1.1}. The dihedral homology $HD(X)$ of any
$D\F$-module $(X,t,r,d,\widetilde{\p})$ over a field of
characteristic zero is isomorphic to the homology of the chain
complex $(M(\overline{X}),b)$.

{\bf Proof}. First, given any $D\F$-module
$(X,t,r,d,\widetilde{\p})$, we note that the specified above
action of the group $\mathbb{Z}_2$ on the chain bicomplex
$(C(\overline{X}),\delta_1,\delta_2)$ induces the action of the
group $\mathbb{Z}_2$ on the chain complex $({\rm
Tot}(C(\overline{X})),\delta)$. This action at any element
$(x_0,\dots,x_n)\in {\rm Tot}(C(\overline{X}))_n$ is given by
$\vartheta(x_0,\dots,x_n)=(\vartheta(x_0),\dots,\vartheta(x_n))$,
where $x_i\in C(\overline{X})_{i,n-i}$ for $0\leqslant i\leqslant
n$. It is easy see that $HD(X)=H(\mathbb{Z}_2;({\rm
Tot}(C(\overline{X})),\delta))$. Now, given any $D\F$-module
$(X,t,r,d,\widetilde{\p})$, consider chain complexes
$(\overline{X},b)$ and $(\overline{X},-b^{'})$. It follows from
the relations $b(1-\overline{T}_n)=(1-\overline{T}_{n-1})b^{'}$,
$n\geqslant 0$, that the chain complex $(L(\overline{X}),b)$,
where $L(\overline{X})_n=\overline{X}_n/{\rm
Im}(1-\overline{T}_n)$ for $n\geqslant 0$, is well defined. The
first and third equalities in $(1.10)$ imply that on the chain
complex $(L(\overline{X}),b)$ correctly acting the group
$\mathbb{Z}_2$ by means of the automorphism
$\vartheta:L(\overline{X})_\bu\to L(\overline{X})_\bu$ of the
order two, which at any element $x\in L(\overline{X})_n$,
$n\geqslant 0$, is given by $\vartheta(x)=\overline{R}_n(x)$. Now,
we consider the chain bicomplex
$$(P'(\overline{X}),d_0,d_1)=({\cal
S}(\mathbb{Z}_2)\otimes_{K[\mathbb{Z}_2]}{\rm
Tot}(C(\overline{X})),d_0,d_1),\quad d_0=1\otimes\delta,$$
$$d_1(a\otimes x)=(-1)^md(a)\otimes x,\quad a\otimes x\in{\cal
S}(\mathbb{Z}_2)_n\otimes_{K[\mathbb{Z}_2]}{\rm
Tot}(C(\overline{X}))_m,\quad n\geqslant 0,\quad m\geqslant 0,$$
where $({\cal S}(\mathbb{Z}_2),d)$ is the specified above standard
free resolvente. Moreover, consider the chain bicomplex
$$(P''(\overline{X}),d_0,d_1)=({\cal
S}(\mathbb{Z}_2)\otimes_{K[\mathbb{Z}_2]}
L(\overline{X}),d_0,d_1),\quad d_0=1\otimes b,$$ $$d_1(a\otimes
x)=(-1)^md(a)\otimes x,\quad a\otimes x\in{\cal
S}(\mathbb{Z}_2)_n\otimes_{K[\mathbb{Z}_2]}
L(\overline{X})_m,\quad n\geqslant 0,\quad m\geqslant 0,$$ where
$({\cal S}(\mathbb{Z}_2),d)$ is the same as above. It was shown in
\cite{Lap8} that over a field of characteristic zero the map of
differential modules $$g:({\rm Tot}(C(\overline{X})),\delta)\to
(L(\overline{X}),b),\quad g(x_0,\dots,x_n)=[x_n]=x_n+{\rm
Im}(1-\overline{T}_n),$$ where $(x_0,\dots,x_n)\in {\rm
Tot}(C(\overline{X}))_n$, induces the isomorphism of homology
modules. Clearly that the chain map $g$ is an
$\mathbb{Z}_2$-equivariant map and, consequently, defines the map
of chain bicomplexes $$G=1\otimes g:(P'(\overline{X}),d_0,d_1)\to
(P''(\overline{X}),d_0,d_1).$$ Now, we consider the spectral
sequences $\{(E_i',d_i)\}_{i\geqslant 0}$ and
$\{(E''_i,d_i)\}_{i\geqslant 0}$ that corresponds to the
bicomplexes $(P'(\overline{X}),d_0,d_1)$ and
$(P''(\overline{X}),d_0,d_1)$, where
$(E'_0,d_0)=(P'(\overline{X}),d_0)$ and
$(E''_0,d_0)=(P''(\overline{X}),d_0)$. Since $G$ is a map of chain
bicomplexes, this map induces the map of spectral sequences
$\{G_i:(E_i',d_i)\to (E_i'',d_i)\}_{i\geqslant 0}$, where
$G_0=1\otimes g$. The map $g$ over a field characteristic zero
induces the isomorphism of homology modules. Therefore the map
$G_0=1\otimes g$ induces the isomorphism of bigraded homology
modules. It implies that the map $G_1:E_1'\to E_1''$ is a
isomorphism of bigraded modules. By using the comparison theorem
of spectral sequences we obtain that all maps $G_i:E_i'\to E_i''$,
$1\leqslant i\leqslant\infty$, are isomorphisms of bigraded
modules. Now, we note that over any field the limit term
$E'_\infty$ of the spectral sequence $\{(E_i',d_i)\}_{i\geqslant
0}$ is isomorphic to the homology of the chain complex that
associated with the chain bicomplex $(P'(\overline{X}),d_0,d_1)$.
Similarly we have that over any field the limit term $E''_\infty$
of the spectral sequence $\{(E_i'',d_i)\}_{i\geqslant 0}$ is
isomorphic to the homology of the chain complex that associated
with the chain bicomplex $(P''(\overline{X}),d_0,d_1)$. It follows
that the dihedral homology $HD(X)$ and the homology of the chain
complex that associated with the chain bicomplex
$(P''(\overline{X}),d_0,d_1)$ are isomorphic. Thus we have the
isomorphism $$HD(X)=H(\mathbb{Z}_2;(L(\overline{X}),b)).$$

Let us show that there is the isomorphism
$H(\mathbb{Z}_2;(L(\overline{X}),b))=H(M(\overline{X}),b)$. Since
we have the isomorphism of modules $P''(\overline{X})_{n.m}={\cal
S}(\mathbb{Z}_2)_n\otimes_{K[\mathbb{Z}_2]}
L(\overline{X})_m=L(\overline{X})_m$, $n\geqslant 0$, $m\geqslant
0$, the chain bicomplex $(P''(\overline{X}),d_0,d_1)$ can be
identified with the chain bicomplex
$(N(\overline{X}),d_0^{\,'},d_1^{\,'})$ that is given by $$
N(\overline{X})_{n,m}=L(\overline{X})_n,\quad n\geqslant 0,\quad
m\geqslant 0,$$
$$d_0^{\,'}=1+(-1)^m\overline{R}_n:N(\overline{X})_{n,m}\to
N(\overline{X})_{n,m-1},\quad d_1^{\,'}=b:N(\overline{X})_{n,m}\to
N(\overline{X})_{n-1,m}.$$ Consider the spectral sequence
$\{(E_i,d_i)\}$ of the chain bicomplex
$(N(\overline{X}),d_0^{\,'},d_1^{\,'})$, where
$(E_0,d_0)=(N(\overline{X}),d_0^{\,'})$. Let us calculate
$(E_1)_{n,m}=H_m(N(\overline{X})_{n,\bu},d_0^{\,'})$. Clearly, we
have $H_0(N(\overline{X})_{n,\bu},d_0^{\,'})=M(\overline{X})_n$.
Let us show that $H_m(N(\overline{X})_{n,\bu},d_0^{\,'})=0$ for
all $m\geqslant 1$. Since we have the equality of chain complexes
$(N(\overline{X})_{n,\bu},d_0^{\,'})=(N(\overline{X})_{n,\bu+1},-d_0^{\,'})$,
it suffices to show that
$$H_1(N(\overline{X})_{n,\bu},d_0^{\,'})={\rm
Ker}(1-\overline{R}_n)/{\rm Im}(1+\overline{R}_n)=0.$$ Suppose
that an element $a\in N(\overline{X})_{n,1}$ satisfies the
condition $d_0^{\,'}(a)=(1-\overline{R}_n)(a)=0$. Consider the
element $c=(1/2)a\in N(\overline{X})_{n,2}=L(\overline{X})_n$.
Since $\overline{R}_n(a)=a$, we have the equality
$d_0^{\,'}(c)=(1/2)(1+\overline{R}_n)(a)=a$, which means that
$H_1(N(\overline{X})_{n,\bu},d_0^{\,'})=0$. Thus, we have
$(E_1)_{n,0}=M(\overline{X})_n$, $n\geqslant 0$, and
$(E_1)_{n,m}=0$, $n\geqslant 0$, $m\geqslant 1$. It is clear that
the differential  $d_1:(E_1)_{n,0}\to (E_1)_{n-1,0}$  induced by
the differential $d_1^{\,'}:N(\overline{X})_{n,0}\to
N(\overline{X})_{n-1,0}$ coincides with the differential
$b:M(\overline{X})_n\to M(\overline{X})_{n-1}$. This, together
with the fact that the limit term
$E_{\infty}=E_2=H(M(\overline{X}),b)$ of the spectral sequence
$\{(E_i,d_i)\}$ over a field is isomorphic to the homology of the
chain complex that associated with the chain bicomplex
$(N(\overline{X}),d_0^{\,'},d_1^{\,'})$, implies that there is the
isomorphism
$H(\mathbb{Z}_2;(L(\overline{X}),b))=H(M(\overline{X}),b)$. Thus,
$HD(X)=H(\mathbb{Z}_2;(L(\overline{X}),b))=H(M(\overline{X}),b)$
and, consequently, we obtain the required isomorphism
$HD(X)=H(M(\overline{X}),b)$.~~~$\blacksquare$

Now, we proceed to the notion of a reflexive module with
$\infty$-simplicial faces.

By a reflexive bigraded module $(X,r)$ we mean any bigraded module
$X$ together with a family of module maps $r=\{r_n:X_{n,\bu}\to
X_{n,\bu}\}$, $n\geqslant 0$, satisfying the conditions
$$r_n^2=1_{X_{n,\bu}},\quad n\geqslant 0.$$ In other words, on
each graded module $X_{n,\bu}$, $n\geqslant 0$, the group
$\mathbb{Z}_2$ of order $2$ with generators $r_n$ acts on the
left.

In what follows, we use the term reflexive differential module for
any triple $(X,r,d)$, where $(X,r)$ is a reflexive bigraded
module, $(X,d)$ is a differential module, and the conditions
$dr_n=r_nd$ for all $n\geqslant 0$ holds.

Now, let us recall that a reflexive module with simplicial faces
\cite{KLS} is defined as a reflexive differential module $(X,r,d)$
together with a family of module maps $\p_i:X_{n,\bu}\to
X_{n-1,\bu}$, $0\leqslant i\leqslant n$, with respect to which the
triple $(X,d,\p_i)$ is a differential module with simplicial faces
and, moreover, the relations $$\p_ir_n=r_{n-1}\p_{n-i},\quad
0\leqslant i\leqslant n,$$ for each $n\geqslant 0$ are true.

It is easy to see that if in the definition of a dihedral module
with simplicial faces we remove the family of automorphisms
$t_n:X_{n,\bu}\to X_{n,\bu}$, $n\geqslant 0$, and the relations
$\p_it_n=t_{n-1}\p_{i-1}$, $0<i\leqslant n$, $\p_0t_n=\p_n$,
$n\geqslant 0$, then we obtain the definition of a reflexive
module with simplicial faces.

{\bf Definition 1.3}. A reflexive module with $\infty$-simplicial
faces or, more briefly, a $R\F$-module, is any quadruple
$(X,r,d,\widetilde{\p})$, where $(X,r,d)$ is a reflexive
differential module and $(X,d,\widetilde{\p})$ is a differential
module with $\infty$-simplicial faces related by the relations
$(1.3)$.

In what follows, we refer to the family of maps
$\widetilde{\p}=\{\p_{(i_1,\dots ,i_k)}:X_{n,\bu}\to
X_{n-k,\bu+k-1}\}$ as the $\F$-dif\-fe\-ren\-tial of the
$R\F$-mo\-du\-le $(X,r,d,\widetilde{\p})$. The maps
$\p_{(i_1,\dots ,i_k)}$ that form the $\F$-dif\-fe\-ren\-tial of a
$R\F$-mo\-du\-le $(X,r,d,\widetilde{\p})$ are called the
$\infty$-simplicial faces of this $R\F$-module.

It is obvious that if in the definition 1.1 we remove the family
of automorphisms $t_n:X_{n,\bu}\to X_{n,\bu}$, $n\geqslant 0$, and
the relations $(1.2)$ then we obtain the definition 1.3. Therefore
each $D\F$-module $(X,t,r,d,\widetilde{\p})$ always defines the
$R\F$-module $(X,r,d,\widetilde{\p})$.

Simple examples of $R\F$-modules are reflexive modules with
simplicial faces. Indeed, given any reflexive module with
simplicial faces $(X,r,d,\p_i)$, we can define an
$\F$-dif\-ferential $\widetilde{\p}=\{\p_{(i_1,\dots ,i_k)}\}$ by
setting $\p_{(i)}=\p_i$, $i\geqslant 0$, and
$\p_{(i_1,\dots,i_k)}=0$, $k>1$, thus obtaining the $R\F$-module
$(X,r,d,\widetilde{\p})$.

Given any reflexive module with $\infty$-simplicial faces
$(X,r,d,\widetilde{\p})$, consider the specified above $\D$-module
$(X,d_0^{\,i})$ and the family of maps
$$R_n=(-1)^{n(n+1)/2}r_n:X_{n,\bu}\to X_{n,\bu},\quad n\geqslant
0.$$ In the same way as it was made in the case of $D\F$-modules
we obtain the relations $d_0^{\,i}R_n=R_{n-i}d_0^{\,i}$,
$i\geqslant 0$, $n\geqslant 0$. Consider the chain complex
$(\overline{X},b)$, which corresponds to the $\D$-module
$(X,d_0^{\,i})$, and the family of maps
$$\overline{R}_n=(R_0+R_1+\dots+R_n):\overline{X}_n\to\overline{X}_n,\quad
n\geqslant 0.$$ By using $d_0^{\,i}R_n=R_{n-i}d_0^{\,i}$,
$i\geqslant 0$, $n\geqslant 0$, we obtain the relations
$b\overline{R}_n=\overline{R}_{n-1}b$, $n\geqslant 0$. This
relations say us that there is a left action of the group
$\mathbb{Z}_2=\{1,\vartheta\,|\,\vartheta^2=1\}$ on the chain
complex $(\overline{X},b)$ of any $R\F$-module
$(X,r,d,\widetilde{\p})$. This left action is defined by means of
the automorphism $\vartheta:\overline{X}\to \overline{X}$ of the
order two, which at any element $x\in \overline{X}_n$ is given by
$\vartheta(x)=\overline{R}_n(x)$.

{\bf Definition 1.4}. We define the reflexive homology $HR(X)$ of
a reflexive module with $\infty$-simplicial faces
$(X,r,d,\widetilde{\p})$ as the hyperhomology
$H(\mathbb{Z}_2;(\overline{X},b))$ of the group $\mathbb{Z}_2$
with coefficients in $(\overline{X},b)$ relative to the specified
above action of the group $\mathbb{Z}_2$ on the chain complex
$(\overline{X},b)$.

It is easy to see that if calculate the hyperhomology
$H(\mathbb{Z}_2;(\overline{X},b))$ by using the specified above
standard free resolvente $({\cal S}(\mathbb{Z}_2),d)$, then we
obtain that the reflexive homology
$HR(X)=H(\mathbb{Z}_2;(\overline{X},b))$ of any $R\F$-module
$(X,r,d,\widetilde{\p})$ is the homology of the chain complex
associated with the chain bicomplex $(R(\overline{X}),b,D)$ that
is defined by $$R(\overline{X})_{n,m}=\overline{X}_n,\quad
n\geqslant 0,\quad m\geqslant 0,$$
$$D=(-1)^n(1+(-1)^m\overline{R}_n):R(\overline{X})_{n,m}\to
R(\overline{X})_{n,m-1},$$
$$b:R(\overline{X})_{n,m}=\overline{X}_n\to \overline{X}_{n-1}=
R(\overline{X})_{n-1,m}.$$ The chain complex associated with the
chain bicomplex $(R(\overline{X}),b,D)$ we denote by $({\rm
Tot}(R(\overline{X})),\widehat{D})$, where $\widehat{D}=b+D$.

Note that if a $R\F$-module $(X,r,d,\widetilde{\p})$ is a
reflexive module with simplicial faces $(X,r,d,\p_i)$, then the
chain bicomplex $(R(\overline{X}),b,D)$ coincides with well-known
\cite{KLS} the chain bicomplex, which computes the reflexive
homology of the reflexive differential module with simplicial
faces $(X,r,d,\p_i)$.

Consider the above chain complex $(\overline{X},b)$ that is
determined by any $R\F$-module $(X,r,d,\widetilde{\p})$. It
follows from the mentioned above relations
$b\overline{R}_n=\overline{R}_{n-1}b$, $n\geqslant 0$, that chain
complex $(N(\overline{X}),b)$, where
$$N(\overline{X})_n=\overline{X}_n\left/{\rm
Im}(1-\overline{R}_n),\quad n\geqslant 0,\right.$$ is well
defined.

The following assertion describes the convenient method of a
computation of the reflexive homology of reflexive modules with
$\infty$-simplicial faces over fields of characteristic zero. This
assertion proved similarly to the theorem 1.1.

{\bf Theorem 1.2}. The reflexive homology $HR(X)$ of any
$R\F$-module $(X,r,d,\widetilde{\p})$ over a field of
characteristic zero is isomorphic to the homology of the chain
complex $(N(\overline{X}),b)$.~~~$\blacksquare$ \vspace{0.5cm}

\centerline{\bf \S\,2. Dihedral and reflexive homology of
involutive $A_\infty$-algebras}
\vspace{0.5cm}

First, following \cite{Kad} (see also \cite{S}), we recall that an
$\A$-algebra $(A,d,\pi_n)$  is any differential module $(A,d)$
with $A=\{A_n\}$, $n\in\mathbb{Z}$, $n\geqslant 0$, $d:A_\bu\to
A_{\bu-1}$, equipped with a family of maps $\{\pi_n:(A^{\otimes
(n+2)})_\bu\to A_{\bu+n}\}$, $n\geqslant 0$, satisfying the
following relations for any integer $n\geqslant 1$:
$$d(\pi_{n-1})=\sum\limits_{m=1}^{n-1}\sum_{t=1}^{m+1}(-1)^{t(n-m)+n+1}
\pi_{m-1}(\underbrace{1\otimes\dots\otimes
1}_{t-1}\otimes\,\pi_{n-m-1}\otimes\underbrace{1\otimes\dots\otimes
1}_{m-t+1}),\eqno(2.1)$$ where
$d(\pi_{n-1})=d\pi_{n-1}+(-1)^n\pi_{n-1}d$. For example, at
$n=1,2,3$ the relations $(2.1)$ take the forms $$d(\pi_0)=0,\quad
d(\pi_1)=\pi_0(\pi_0\otimes 1)-\pi_0(1\otimes\pi_0),$$
$$d(\pi_2)=\pi_0(\pi_1\otimes 1+1\otimes\pi_1)-\pi_1(\pi_0\otimes
1\otimes 1-1\otimes\pi_0\otimes 1+1\otimes 1\otimes\pi_0).$$

Now, we recall the notion of an involutive $\A$-algebra \cite{B},
which generalizes the notion of a differential associative algebra
with a involution. An involutive $\A$-al\-geb\-ra $(A,d,\pi_n,*)$
is defined as an $\A$-algebra $(A,d,\pi_n)$ together with the
automorphism of graded modules $*:A_\bu\to A_\bu$ (the notation
$*(a)=a^*$ for $a\in A$), which at any elements $a\in A$ and
$a_0,a_1,\dots,a_n,a_{n+1}\in A$ satisfies the conditions
$$a^{**}=a,\quad d(a^*)=d(a)^*,$$ $$\pi_n(a_0\otimes
a_1\otimes\dots\otimes a_n\otimes
a_{n+1})^*=(-1)^\varepsilon\pi_n(a_{n+1}^*\otimes
a_n^*\otimes\dots\otimes a_1^*\otimes a_0^*),\eqno(2.2)$$
$$\varepsilon=\frac{n(n+1)}{2}+\sum_{0\leqslant i<j\leqslant
n+1}|a_i||a_j|,\quad n\geqslant 0,$$ where  $|a|=q$ means that
$a\in A_q$.

Given any involutive $\A$-algebra $(A,d,\pi_n,*)$, consider the
dihedral differential module $(^{\varrho}{\cal M}(A),t,r,d)$,
where $\varrho=\pm 1$, defined by $$^{\varrho}{\cal
M}(A)=\{^{\varrho}{\cal M}(A)_{n,m}\},\quad ^{\varrho}{\cal
M}(A)_{n,m}=(A^{\otimes (n+1)})_m,\quad n\geqslant 0,\quad
m\geqslant 0,$$ $$t_n(a_0\otimes\dots\otimes
a_n)=(-1)^{|a_n|(|a_0|+\dots+|a_{n-1}|)}a_n\otimes a_0\otimes
a_1\otimes \dots\otimes a_{n-1},$$ $$r_n(a_0\otimes\dots\otimes
a_n)=\varrho(-1)^{\sum_{0<i<j\leqslant n}|a_i||a_j|}a_0^*\otimes
a_n^*\otimes a_{n-1}^*\dots\otimes a_1^*,$$
$$d(a_0\otimes\dots\otimes
a_n)=\sum_{i=0}^n(-1)^{|a_0|+\dots+|a_{i-1}|}a_0\otimes\dots\otimes
a_{i-1}\otimes d(a_i)\otimes a_{i+1}\otimes\dots\otimes a_n.$$ It
is easy to verify that maps $t_n$ and $r_n$ satisfy the conditions
$t_n^{n+1}=1$, $r_n^2=1$, $t_nr_n=r_nt_n^{-1}$ for all $n\geqslant
0$. Therefore, $(^{\varrho}{\cal M}(A),t,r,d)$ is indeed a
dihedral differential module. Now, consider the family of maps
$\widetilde{\p}=\{\p_{(i_1,\dots,i_k)}:{^{\varrho}{\cal
M}(A)_{n,p}}\to {^{\varrho}{\cal M}(A)_{n-k,p+k-1}}\}$,
$0\leqslant i_1<\dots<i_k\leqslant n$, $n\geqslant 0$, $p\geqslant
0$, defined by $$\p_{(i_1,\dots,i_k)}=$$
$$=\left\{\begin{array}{ll} (-1)^{k(p-1)}1^{\otimes
j}\otimes\,\pi_{k-1}\otimes 1^{\otimes
(n-k-j)}\,,\,\,\,\mbox{if}\,\,0\leqslant j\leqslant
n-k&\\\quad\mbox{and}\,\,\, (i_1,\dots,i_k)=(j,j+1,\dots,j+k-1);\\
(-1)^{q(k-1)}\p_{(0,1,\dots,k-1)}t_n^q\,,\,\,\,\mbox{if}\,\,1\leqslant
q\leqslant k&\\\quad\mbox{and}\,\,\,
(i_1,\dots,i_k)=(0,1,\dots,k-q-1,n-q+1,n-q+2,\dots,n);\\
0,\quad\mbox{otherwise}.&\\
\end{array}\right.\eqno(2.3)$$

{\bf Theorem 2.1}. For any involutive $\A$-algebra
$(A,d,\pi_n,*)$, the specified above five-tuple $(^{\varrho}{\cal
M}(A),t,r,d,\widetilde{\p})$ is a dihedral module with
$\infty$-simplicial faces.

{\bf Proof}. In \cite{Lap8} it was shown that the quadruple
$(^{\varrho}{\cal M}(A),t,d,\widetilde{\p})$ is a cyclic
differential module with $\infty$-simplicial faces. Therefore, it
is remains only to check that the maps
$\p_{(i_1,\dots,i_k)}\in\widetilde{\p}$ defined by $(2.3)$ satisfy
the relations $(1.3)$. Consider several cases.

1). Suppose that $(i_1,\dots,i_k)=(j,j+1,\dots,j+k-1)$, where
$1\leqslant j\leqslant n-k$. In this case, on the one hand, we
have at any element $a_0\otimes\dots\otimes a_n\in (A^{n+1})_p$
the equalities
$$\p_{(j,j+1,\dots,j+k-1)}r_n(a_0\otimes\dots\otimes
a_n)=\varrho(-1)^\alpha\p_{(j,j+1,\dots,j+k-1)}(a_0^*\otimes
a_n^*\otimes\dots\otimes a_1^*)=$$
$$=\varrho(-1)^{\alpha+k(p-1)}(1^{\otimes
j}\otimes\pi_{k-1}\otimes 1^{\otimes (n-k-j)})(a_0^*\otimes
a_n^*\otimes\dots\otimes a_1^*)=$$
$$=\varrho(-1)^{\alpha+k(p-1)+\beta}a_0^*\otimes
a_n^*\otimes\dots\otimes a_{n-j+2}^*\,\otimes$$
$$\otimes\,\pi_{k-1}(a_{n-j+1}^*\otimes\dots\otimes
a_{n-j-k+1}^*)\otimes a_{n-j-k}^*\otimes\dots\otimes a_1^*,$$
where $$\alpha=\sum\limits_{0<i<j\leqslant n}|a_i||a_j|,\quad
\beta=(k-1)(|a_0|+|a_n|+\dots+|a_{n-j+2}|).$$ On the other hand,
we have the equalities
$$(-1)^{k(k-1)/2}r_{n-k}\p_{(n-j-k+1,n-j-k+2,\dots,n-j)}(a_0\otimes\dots\otimes
a_n)=$$ $$=(-1)^{(k(k-1)/2)+k(p-1)}r_{n-k}(1^{\otimes
(n-j-k+1)}\otimes\pi_{k-1}\otimes 1^{\otimes
(j-1)})(a_0\otimes\dots\otimes a_n)=$$
$$=(-1)^{(k(k-1)/2)+k(p-1)+\gamma}r_{n-k}(a_0\otimes\dots\otimes
a_{n-j-k}\otimes\pi_{k-1}(a_{n-j-k+1}\otimes\dots\otimes
a_{n-j+1})\,\otimes$$ $$\otimes\,a_{n-j+2}\otimes\dots\otimes
a_n)=$$
$$=\varrho(-1)^{(k(k-1)/2)+k(p-1)+\gamma+\delta}a_0^*\otimes
a_n^*\otimes\dots\otimes a_{n-j+2}^*\otimes$$
$$\otimes\,\pi_{k-1}(a_{n-j-k+1}\otimes\dots\otimes
a_{n-j+1})^*\otimes a_{n-j-k}^*\otimes\dots\otimes a_1^*=$$
$$=\varrho(-1)^{(k(k-1)/2)+k(p-1)+\gamma+\delta+\mu}a_0^*\otimes
a_n^*\otimes\dots\otimes a_{n-j+2}^*\,\otimes$$
$$\otimes\,\pi_{k-1}(a_{n-j+1}^*\otimes\dots\otimes
a_{n-j-k+1}^*)\otimes a_{n-j-k}^*\otimes\dots\otimes a_1^*,$$
where $$\gamma=(k-1)(|a_0|+|a_1|+\dots+|a_{n-j-k}|),$$
$$\delta=\sum\limits_{0<i<j\leqslant
n}|a_i||a_j|-\sum\limits_{n-j-k+1\leqslant s<t\leqslant
n-j+1}|a_s||a_t|\,+$$
$$+\,(k-1)(|a_1|+\dots+|a_{n-j-k}|+|a_{n-j+2}|+\dots+|a_n|),$$
$$\mu=\frac{k(k-1)}{2}+\sum\limits_{n-j-k+1\leqslant s<t\leqslant
n-j+1}|a_s||a_t|.$$ Since $\alpha+\beta\equiv
(k(k-1)/2)+\gamma+\delta+\mu\,{\rm mod}(2)$, we obtain we obtain
the required relation
$$\p_{(j,j+1,\dots,j+k-1)}r_n=(-1)^{k(k-1)/2}r_{n-k}\p_{(n-j-k+1,n-j-k+2,\dots,n-j)}.$$

2). Suppose that $(i_1,\dots,i_k)=(0,1,\dots,k-1)$. In this case,
on the one hand, we have at any element $a_0\otimes\dots\otimes
a_n\in (A^{n+1})_p$  the equalities
$$\p_{(0,1,\dots,k-1)}r_n(a_0\otimes\dots\otimes
a_n)=\varrho(-1)^\alpha\p_{(0,1,\dots,k-1)}(a_0^*\otimes
a_n^*\otimes\dots\otimes a_1^*)=$$
$$=\varrho(-1)^{\alpha+k(p-1)}(\pi_{k-1}\otimes 1^{\otimes
(n-k)})(a_0^*\otimes a_n^*\otimes\dots\otimes a_1^*)=$$
$$\varrho(-1)^{\alpha+k(p-1)}\pi_{k-1}(a_0^*\otimes
a_n^*\otimes\dots\otimes a_{n-k+1}^*)\otimes
a_{n-k}^*\otimes\dots\otimes a_1^*,$$ where $\alpha$ is the same
as in 1). On the other hand, since the second equality in $(2.3)$
at $q=k$  becomes the equalities
$\p_{(n-k+1,n-k+2,\dots,n)}=\p_{(0,1,\dots,k-1)}t_n^k$, we have
the equality
$$(-1)^{k(k-1)/2}r_{n-k}\p_{(n-k+1,n-k+2,\dots,n)}(a_0\otimes\dots\otimes
a_n)=$$
$$=(-1)^{k(k-1)/2}r_{n-k}\p_{(0,1,\dots,k-1)}t_n^k(a_0\otimes\dots\otimes
a_n)=$$
$$=(-1)^{(k(k-1)/2)+\nu}r_{n-k}\p_{(0,1,\dots,k-1)}(a_{n-k+1}\otimes\dots\otimes
a_n\otimes a_0\otimes\dots\otimes a_{n-k})=$$
$$=(-1)^{(k(k-1)/2)+\nu+k(p-1)}r_{n-k}(\pi_{k-1}\otimes 1^{\otimes
(n-k)})(a_{n-k+1}\otimes\dots\otimes a_n\otimes
a_0\otimes\dots\otimes a_{n-k})=$$
$$=(-1)^{(k(k-1)/2)+\nu+k(p-1)}r_{n-k}(\pi_{k-1}(a_{n-k+1}\otimes\dots\otimes
a_n\otimes a_0)\otimes a_1\otimes\dots\otimes a_{n-k})=$$
$$=\varrho(-1)^{(k(k-1)/2)+\nu+k(p-1)+\vartheta}\pi_{k-1}(a_{n-k+1}\otimes\dots\otimes
a_n\otimes a_0)^*\otimes a_{n-k}^*\otimes\dots\otimes a_1^*=$$
$$=\varrho(-1)^{(k(k-1)/2)+\nu+k(p-1)+\vartheta+\lambda}\pi_{k-1}(a_0^*\otimes
a_n^*\otimes a_{n-k+1}^*)\otimes a_{n-k}^*\otimes\dots\otimes
a_1^*,$$ where
$$\nu=(|a_{n-k+1}|+\dots+|a_n|)(|a_0|+\dots+|a_{n-k}|),\quad
\vartheta=\sum_{0<i<j\leqslant n-k}|a_i||a_j|,$$
$$\lambda=\frac{k(k-1)}{2}+\sum_{n-k+1\leqslant i<j\leqslant
n}|a_i||a_j|+(|a_{n-k+1}|+\dots+|a_n|)|a_0|.$$ Since
$(k(k-1)/2)+\nu+\vartheta+\lambda\equiv\alpha\,{\rm mod}(2)$, we
obtain the required relation
$$\p_{(0,1,\dots,k-1)}r_n=(-1)^{k(k-1)/2}r_{n-k}\p_{(n-k+1,n-k+2,\dots,n)}.$$

3). Suppose that
$(i_1,\dots,i_k)=(0,1,\dots,k-q-1,n-q+1,n-q+2,\dots,n)$, where
$1\leqslant q\leqslant k$. The second equality in $(2.3)$ follows
that
$$\p_{(0,1,\dots,k-q-1,n-q+1,n-q+2,\dots,n)}=(-1)^{q(k-1)}\p_{(0,1,\dots,k-1)}t_n^q.$$
By using
$t_n^qr_n=r_nt_n^{-q}=r_nt_n^{n+1-q}=r_nt_n^{n-k+1}t_n^{k-q}$, we
have
$$\p_{(0,1,\dots,k-q-1,n-q+1,n-q+2,\dots,n)}r_n=(-1)^{q(k-1)}\p_{(0,1,\dots,k-1)}t_n^qr_n=$$
$$=(-1)^{q(k-1)}\p_{(0,1,\dots,k-1)}r_nt_n^{n-k+1}t_n^{k-q}=$$ $$=
(-1)^{q(k-1)+(k(k-1)/2)}r_{n-k}\p_{(n-k+1,\dots,n)}t_n^{n-k+1}t_n^{k-q}=$$
$$=(-1)^{q(k-1)+(k(k-1)/2)}r_{n-k}t_{n-k}^{n-k+1}\p_{(0,1,\dots,k-1)}t_n^{k-q}=$$
$$=(-1)^{q(k-1)+(k(k-1)/2)}r_{n-k}\p_{(0,1,\dots,k-1)}t_n^{k-q}.$$
The second equality in $(2.3)$ follows that
$$\p_{(0,1,\dots,k-1)}t_n^{k-q}=(-1)^{(k-q)(k-1)}\p_{(0,1,\dots,q-1,n-k+q+1,\dots,n)}.$$
Since $q(k-1)+(k(k-1)/2)+(k-q)(k-1)\equiv (k(k-1)/2)\,{\rm
mod}(2)$, we obtain the required relation
$$\p_{(0,1,\dots,k-q-1,n-q+1,n-q+2,\dots,n)}r_n=
(-1)^{k(k-1)/2}r_{n-k}\p_{(0,1,\dots,q-1,n-k+q+1,\dots,n)}.$$
Thus, all maps $\p_{(i_1,\dots,i_k)}\in\widetilde{\p}$ satisfy the
conditions $(1.3)$ and, consequently, the five-tuple
$(^{\varrho}{\cal M}(A),t,r,d,\widetilde{\p})$ is a
$D\F$-module.~~~$\blacksquare$

Note that if an involutive $\A$-algebra $(A,d,\pi_n,*)$ is a
differential involutive associative algebra $(A,d,\pi,*)$, where
$\pi=\pi_0$ and $\pi_n=0$, $n>0$, then the $D\F$-module
$(^{\varrho}{\cal M}(A),t,r,d,\widetilde{\p})$ coincides with the
well-known \cite{KLS} dihedral module with simplicial faces that
defined by a differential involutive associative algebra
$(A,d,\pi,*)$.

{\bf Definition 2.1}. We define the dihedral homology
$^{\varrho}HD(A)$ of an arbitrary involutive $\A$-al\-gebra
$(A,d,\pi_n,*)$ as the dihedral homology $HD(^\varrho{\cal M}(A))$
of the $D\F$-mo\-dule $(^\varrho{\cal
M}(A),t,r,d,\widetilde{\p})$.

Thus, the dihedral homology $^{\varrho}HD(A)$ of an involutive
$\A$-algebra $(A,d,\pi_n,*)$ is the homology of the chain complex
$({\rm Tot}(D(\overline{^{\varrho}{\cal
M}(A)})),\widehat{\delta})$ associated with the triple chain
complex $(D(\overline{^{\varrho}{\cal
M}(A)}),\delta_1,\delta_2,\delta_3)$.

Note that if an involutive $\A$-algebra $(A,d,\pi_n,*)$ is a
differential involutive associative algebra $(A,d,\pi,*)$, where
$\pi=\pi_0$ and $\pi_n=0$, $n>0$, then the triple chain complex
$(D(\overline{^{\varrho}{\cal M}(A)}),\delta_1,\delta_2,\delta_3)$
coincides with the well-known \cite{KLS} triple chain complex that
computes the dihedral homology of a differential involutive
associative algebra $(A,d,\pi,*)$.

If apply the theorem 1.1 to the $D\F$-module $(^{\varrho}{\cal
M}(A),t,r,d,\widetilde{\p})$, then we obtain the following
assertion.

{\bf Corollary 2.1}. The dihedral homology $^{\varrho}HD(A)$ of
any involutive $\A$-algebra $(A,d,\pi_n,*)$ over a field of
characteristic zero is isomorphic to the homology of the chain
complex $(M(\overline{^{\varrho}{\cal
M}(A)}),b)$.~~~$\blacksquare$

Note that the chain complex $(M(\overline{^{\varrho}{\cal
M}(A)}),b)$ was defined in \cite{B} without employing dihedral
modules with $\infty$-simplicial faces. In \cite{B}, where the
ground ring is a field of characteristic zero, the homology of
this chain complex was referred to as the dihedral homology of the
involutive $\A$-algebra $(A,d,\pi_n,*)$. Thus, the corollary 2.1
implies that, over fields of characteristic zero, the definition
2.1 of the dihedral homology of an involutive $\A$-algebra is
equivalent to that given in \cite{B}.

As it was said above, an arbitrary $D\F$-module always can be
considered as a $R\F$-module. Therefore, the $D\F$-module
$(^{\varrho}{\cal M}(A),t,r,d,\widetilde{\p})$ that defined by any
involutive $\A$-algebra $(A,d,\pi_n,*)$ always determines the
$R\F$-module $(^{\varrho}{\cal M}(A),r,d,\widetilde{\p})$.

Note that if an involutive $\A$-algebra $(A,d,\pi_n,*)$ is a
differential involutive associative algebra $(A,d,\pi,*)$, where
$\pi=\pi_0$ and $\pi_n=0$, $n>0$, then the $R\F$-module
$(^{\varrho}{\cal M}(A),r,d,\widetilde{\p})$ coincides with the
well-known \cite{KLS} reflexive module with simplicial faces that
defined by a differential involutive associative algebra
$(A,d,\pi,*)$.

{\bf Definition 2.1}. We define the reflexive homology
$^{\varrho}HR(A)$ of an arbitrary involutive $\A$-al\-gebra
$(A,d,\pi_n,*)$ as the reflexive homology $HR({\cal M}(A))$ of the
$R\F$-module $(^\varrho{\cal M}(A),r,d,\widetilde{\p})$.

Thus, the reflexive homology $^{\varrho}HR(A)$ of an involutive
$\A$-algebra $(A,d,\pi_n,*)$ is the homology of the chain complex
$({\rm Tot}(R(\overline{^{\varrho}{\cal M}(A)})),\widehat{D})$
associated with the chain bicomplex $(R(\overline{^{\varrho}{\cal
M}(A)}),b,D)$.

Note that if an involutive $\A$-algebra $(A,d,\pi_n,*)$ is a
differential involutive associative algebra $(A,d,\pi,*)$, where
$\pi=\pi_0$ and $\pi_n=0$, $n>0$, then the chain bicomplex
$(R(\overline{^{\varrho}{\cal M}(A)}),b,D)$ coincides with the
well-known \cite{KLS} chain bicomplex that computes the reflexive
homology of a differential involutive associative algebra
$(A,d,\pi,*)$.

If apply the theorem 1.2 to the $R\F$-module $(^{\varrho}{\cal
M}(A),r,d,\widetilde{\p})$, then we obtain the following
assertion.

{\bf Corollary 2.2}. The reflexive homology $^{\varrho}HR(A)$ of
any involutive $\A$-algebra $(A,d,\pi_n,*)$ over a field of
characteristic zero is isomorphic to the homology of the chain
complex $(N(\overline{^{\varrho}{\cal
M}(A)}),b)$.~~~$\blacksquare$
\vspace{0.5cm}

\centerline{\bf \S\,3. An exact sequence for the dihedral and the
reflexive} \centerline{\bf homology of involutive homotopy unital
$\A$-algebras}
\vspace{0.5cm}

To introduce the notion of an involutive homotopy unital
$\A$-algebra, we need definitions and constructions related to the
notions of a (nonsymmetric) operad and an algebra over an operad
in the category of differential modules (see, e.g., \cite{LV}).

A differential family or, more briefly, a family, ${\cal
E}=\{{\cal E}(j) \}_{j\geqslant 0}$ is any family of differential
modules $({\cal E}(j),d)$, $j\geqslant 0$. A morphism of families
$f:{\cal E}'\to {\cal E}''$ is any family of maps
$\alpha=\{\alpha(j):({\cal E}'(j),d)\to ({\cal
E}''(j),d)\}_{j\geqslant 0}$ of differential modules. For any
families ${\cal E}'$ and ${\cal E}''$, the family ${\cal
E}'\times{\cal E}''$ is defined by $$({\cal E}'\times{\cal
E}'')(j)=\bigoplus\limits_{j_1+\cdots +j_k=j}{\cal
E}'(k)\otimes{\cal E}''(j_1) \otimes\cdots\otimes{\cal
E}''(j_k),\quad j\geqslant 0.$$ Clearly, the $\times$-product of
families thus defined is associative i.e., given any families
${\cal E}$, ${\cal E}'$ and ${\cal E}''$, there is an isomorphism
of families ${\cal E}\times({\cal E}'\times{\cal
E}'')\approx({\cal E}\times{\cal E}') \times{\cal E}''$.

A (nonsymmetric) operad $({\cal E},\gamma)$ is any family ${\cal
E}$ together with a morphism of families $\gamma:{\cal E}\times
{\cal E}\to {\cal E}$ satisfying the condition
$\gamma(\gamma\times 1)=\gamma(1\times\gamma)$. Moreover, there
exists an element $1\in{\cal E}(1)_0$ such that $\gamma(1\otimes
e_j)=e_j$ for each $e_j\in {\cal E}(j)$, $j\geqslant 0$, and
$\gamma(e_j\otimes 1\otimes\dots\otimes 1)=e_j$ for each $e_j\in
{\cal E}(j)$, $j\geqslant 1$. In what follows, we write elements
of the form $\gamma(e_k\otimes e_{j_1}\otimes\dots\otimes
e_{j_k})$ as $e_k(e_{j_1}\otimes\dots\otimes e_{j_k})$.

A morphism of operads $f:{\cal E}',\gamma)\to ({\cal E}'',\gamma)$
is a morphism of families $f:{\cal E}'\to {\cal E}''$ satisfying
the condition $f\gamma=\gamma(f\times f)$.

A canonical example of an operad is the operad $({\cal
E}_X,\gamma)$ defined by $$({\cal E}_X(j),d)=(\Hom(X^{\otimes
j};X),d),\quad \gamma(f_k\otimes f_{j_1}\otimes\dots\otimes
f_{j_k})=f_k(f_{j_1}\otimes\dots\otimes f_{j_k})$$ for any
differential module $(X,d)$.

An algebra over on operad $({\cal E},\gamma)$ or, more briefly, an
${\cal E}$-algebra $(X,d,\alpha)$ is defined as a differential
module $(X,d)$ together with a fixed morphism of operads
$\alpha:{\cal E}\to {\cal E}_X$.

An important example of an operad is the Stasheff operad
$(\A,\gamma)$. As a graded operad this is the free operad with
generators $\pi_n\in\A(n+2)_n$, $n\geqslant 0$; at each generator
$\pi_{n-1}$, $n\geqslant 1$, the differential is defined by
$(2.1)$.

It is easy to see that, given a differential module $(A,d)$, where
$A=\{A_n\}$, $n\geqslant 0$, $d:A_\bu\to A_{\bu-1}$, defining the
structure of an algebra $(A,d,\alpha)$ over the operad
$(\A,\gamma)$ on the differential module $(A,d)$ is equivalent to
specifying a family of maps $$\{\pi_n=\alpha(\pi_n):(A^{\otimes
(n+2)})_\bu\to A_{\bu+n}\},\quad n\geqslant 0,$$ for which the
relations $(2.1)$ hold, i.e., to endowing $(A,d)$ with the
structure of the $\A$-algebra $(A,d,\pi_n)$.

Now, following \cite{Lyub} (see also \cite{Lap6}), we recall the
notion of a homotopy unital $\A$-al\-gebra. Consider the operad
$(\A^{su}<\!u,h\!>,\gamma)$ introduced in \cite{Lyub}. As a graded
operad, $(\A^{su}<\!u,h\!>,\gamma)$ has the generators
$$\pi_n\in(\A^{su}<\!u,h\!>)(n+2)_n,\quad n\geqslant 0,\quad
1^{su}\in (\A^{su}<\!u,h\!>)(0)_0,$$ $$u\in
(\A^{su}<\!u,h\!>)(0)_0,\quad h\in (\A^{su}<\!u,h\!>)(0)_1,$$
which satisfy the relations $$\pi_0(1^{su}\otimes 1)=1,\quad
\pi_0(1\otimes 1^{su})=1,\quad\pi_n(1^{\otimes k}\otimes
1^{su}\otimes 1^{\otimes (n-k+1)})=0,\quad n>0,\eqno(3.1)$$ where
$0\leqslant k\leqslant n+1$; the differential is defined at the
generators by the formulae $(2.1)$ and $$d(1^{su})=0,\quad
d(u)=0,\quad d(h)=u-1^{su}.\eqno(3.2)$$

The operad $(\A^{su}<\!u,h\!>,\gamma)$ contains the suboperad
$(\A^{hu},\gamma)$ with generators $$\tau_0^0=u\in
(\A^{su}<\!u,h\!>)(0)_0,\quad
\tau^\varnothing_n=\pi_{n-1}\in(\A^{su}<\!u,h\!>)(n+1)_{n-1},\quad
n\geqslant 1,$$ $$\tau_n^{j_q,\dots,j_1}=$$
$$=\pi_{n-1}(\underbrace{\overbrace{\underbrace{1^{\otimes
n_1}}_{j_1}\otimes h\otimes 1^{\otimes n_2}}^{j_2}\otimes h\otimes
1^{\otimes n_3}\dots\otimes 1^{n_k}\otimes h\otimes 1^{\otimes
n_{k+1}}\otimes\dots\otimes 1^{\otimes n_q}}_{j_q}\otimes h\otimes
1^{\otimes n_{q+1}})\in$$
$$\in(\A^{su}<\!u,h\!>)(n-q+1)_{n+q-1},\quad n\geqslant 1,\quad
q\geqslant 1,\quad n\geqslant j_q>\dots>j_1\geqslant 0,\quad
n_s\geqslant 0,$$ $$1\leqslant s\leqslant q+1,\quad
n_1+\dots+n_{q+1}=n-q+1,\quad j_k=n_1+\dots+n_k+k-1,\quad
1\leqslant k\leqslant q.$$

It is worth mentioning for clarity that each $j_k$, $1\leqslant
k\leqslant q$, is equal to the number of all tensor factors on the
left of the $k$th occurrence of $h$, counting from the beginning
of the tensor stack to the right. For example, according to this
rule, we have $$\tau_4^{3,1}=\pi_3(1\otimes h\otimes 1\otimes
h\otimes 1),\quad \tau_3^{3,2}=\pi_2(1\otimes 1\otimes h\otimes
h),$$ $$\tau_5^{5,1,0}=\pi_4(h\otimes h\otimes 1\otimes 1\otimes
1\otimes h).$$
Note that, in \cite{Lyub}, an element
$\tau_n^{j_q,\dots,j_1}$ is denoted by
$m_{n_1,n_2,\dots,n_{q+1}}$, where the numbers $n_1,\dots,n_{q+1}$
are the same as in the above expression for
$\tau_n^{j_q,\dots,j_1}$. At the generators
$\tau_n^{j_q,\dots,j_1}$, $n\geqslant 0$, $q\geqslant 0$,
$n+q\geqslant 1$, where
$\tau_n^{j_q,\dots,j_1}=\tau^\varnothing_n=\pi_{n-1}$ for $q=0$
and $n\geqslant 1$, the differential is completely determined by
the formulae $(2.1)$, $(3.2)$ and the relations $(3.1)$.

Algebras $(A,d,\alpha)$ over the operad $(\A^{hu},\gamma)$, i.e.,
$\A^{hu}$-algebras, are called homotopy unital $\A$-algebras.

It is easy to see that defining the structure of an
$\A^{hu}$-algebra $(A,d,\alpha)$ on a differential module $(A,d)$,
where $A=\{A_n\}$, $n\in\mathbb{Z}$, $n\geqslant 0$, $d:A_\bu\to
A_{\bu-1}$, is equivalent to specifying a family of maps
$$\{\alpha(\tau_n^{j_q,\dots,j_1}):(A^{\otimes (n-q+1)})_\bu\to
A_{\bu+n+q-1}\,\,|\,\,n\in\mathbb{Z},\,\, n\geqslant
0,\,\,q\geqslant 0,\,\,n+q\geqslant 1\},$$ $$n-q+1\geqslant
0,\quad j_q,\dots,j_1\in\mathbb{Z},\quad n\geqslant
j_q>\dots>j_1\geqslant 0,$$ which satisfy the relations
$$d(\alpha(\tau_n^{j_q,\dots,j_1}))=\alpha(d(\tau_n^{j_q,\dots,j_1})),\eqno(3.3)$$
where
$d(\alpha(\tau_n^{j_q,\dots,j_1}))=d\alpha(\tau_n^{j_q,\dots,j_1})+(-1)^{n+q}\alpha(\tau_n^{j_q,\dots,j_1})d$.
In what follows, we denote the map
$\alpha(\tau_n^{j_q,\dots,j_1})$ by $\tau_n^{j_q,\dots,j_1}$. For
example, by $(3.3)$ we have the following relations:
$$d(\tau_2^0)=\tau_1^0\pi_0+\pi_0(\tau_1^0\otimes
1)-\pi_1(\tau_0^0\otimes 1\otimes 1),$$
$$d(\tau_2^{2,0})=-\tau_1^1\tau_1^0-\tau_1^0\tau_1^1-\tau_2^2(\tau_0^0\otimes
1)+\tau_2^0(1\otimes\tau_0^0),$$
$$d(\tau_3^3)=\tau_1^\varnothing(1\otimes\tau_2^2)-\tau_2^\varnothing(1\otimes
1\otimes\tau_1^1)-\tau_1^1\tau_2^\varnothing-\tau_2^2(\tau_1^\varnothing\otimes
1)\,+$$
$$+\,\tau_2^2(1\otimes\tau_1^\varnothing)+\tau_3^\varnothing(1\otimes
1\otimes 1\otimes\tau_0^0).$$

Note that, for $q=0$, the relations $(3.3)$ in which the maps
$\alpha(\tau_n^\varnothing)$, $n\geqslant 1$, are denoted by
$\pi_{n-1}$ transforms into the relations $(2.1)$. Thus, each
$\A^{hu}$-algebra $(A,d,\tau_n^{j_q,\dots,j_1})$ is an
$\A$-algebra $(A,d,\pi_n)$. In what follows, we write
$\A^{hu}$-algebras $(A,d,\tau_n^{j_q,\dots,j_1})$ in the form
$(A,d,\pi_n,\tau_n^{j_q,\dots,j_1})$, where
$\pi_n=\tau^\varnothing_{n+1}$, $n\geqslant 0$.

Note also that, by virtue of $(3.2)$ and $(3.3)$, we have
$$d(\tau_0^0)=0,\quad d(\tau_1^0)=\pi_0(\tau_0^0\otimes 1)-1,\quad
d(\tau_1^1)=\pi_0(1\otimes\tau_0^0)-1,$$ for any $\A^{hu}$-algebra
$(A,d,\pi_n,\tau_n^{j_q,\dots,j_1})$; thus, the map $\tau_0^0:K\to
A$ or, to be more precise, the element $\tau_0^0(1)\in A_0$ is up
to homotopy the unit of the differential homotopy associative
algebra $(A,d,\pi_0)$. Moreover, it was shown in \cite{Lap8} that
in an explicit form the relations $(3.3)$, for $q=1$ and
$j_q=j_1=n\geqslant 2$, are written as the following equalities:
$$d(\tau_n^n)=\sum\limits_{m=1}^{n-1}\sum_{t=0}^{m-1}(-1)^{t(n-m)+n}
\tau_m^m (\underbrace{1\otimes\dots\otimes
1}_t\otimes\,\pi_{n-m-1}\otimes\underbrace{1\otimes\dots\otimes
1}_{m-t-1})\,+$$ $$+\sum\limits_{m=1}^n(-1)^{mn+1}
\pi_{m-1}(1\otimes\dots\otimes
1\otimes\tau_{n-m}^{n-m}).\eqno(3.4)$$

Now we introduce the notion of an involutive homotopy unital
$\A$-algebra, which generalizes the notion of a differential
unital associative algebra with a involution.

{\bf Definition 3.1}. We define the involutive homotopy unital
$\A$-algebra, briefly, the involutive $\A^{hu}$-algebra
$(A,d,\pi_n,\tau_n^{j_q,\dots,j_1},*)$ as any five-tuple
$(A,d,\pi_n,\tau_n^{j_q,\dots,j_1},*)$ that satisfy the following
conditions:

1). $(A,d,\pi_n,*)$ is an involutive $\A$-algebra;

2). $(A,d,\pi_n,\tau_n^{j_q,\dots,j_1})$ is a homotopy unital
$\A$-algebra;

3). For any elements $a_0,\dots,a_{n-q}\in A$, $n\geqslant 1$,
$q\geqslant 1$, there are the relations
$$\tau_n^{j_q,\dots,j_1}(a_0\otimes a_1\otimes\dots\otimes
a_{n-q-1}\otimes a_{n-q})^*=$$
$$=(-1)^\varepsilon\tau_n^{n-j_1,\dots,n-j_q}(a_{n-q}^*\otimes
a_{n-q-1}^*\otimes\dots\otimes a_1^*\otimes a_0^*),$$ $$n\geqslant
j_q>\dots>j_1\geqslant 0,\quad
\varepsilon=\frac{n(n-1)}{2}+\frac{q(q-1)}{2}+\sum_{0\leqslant
i<j\leqslant n-q}|a_i||a_j|,$$ where  $|a|=q$ means that $a\in
A_q$.

Given an $\A^{hu}$-algebra $(A,d,\pi_n,\tau_n^{j_q,\dots,j_1},*)$
consider the family of maps $$s^{k-1}:(A^{\otimes (n+1)})_p\to
(A^{\otimes (n-k+2)})_{p+k},\quad n\geqslant 0,\quad p\geqslant
0,\quad 0\leqslant k\leqslant n+1,$$ defined by the following
formulae:
$$s^{k-1}=(-1)^\varepsilon\underbrace{1\otimes\dots\otimes
1}_{n-k+1}\otimes\,\tau_k^k,\eqno(3.5)$$
$$\varepsilon=(k-1)p+(n-k+1)k+\frac{k(k-1)}{2}+n+1.$$ Consider
also the $D\F$-module $(^\varrho{\cal M}(A),t,r,d,\widetilde{\p})$
corresponding to this $\A^{hu}$-algebra. As mentioned above, this
$D\F$-module determines the $\D$-module $(^\varrho{\cal
M}(A),d_1^{\,i})$ according to $(1.4)$ with $q =1$. In
\cite{Lap8}, it was shown that, given any integer $k\geqslant 0$,
the families of maps $$\{d_1^{\,i}:{^\varrho{\cal
M}(A)}_{*,\bu}\to {^\varrho{\cal
M}(A)}_{*-i,\bu+i-1}\}\quad\mbox{and}\quad
\{s^{i-1}:{^\varrho{\cal M}(A)}_{*,\bu}\to {^\varrho{\cal
M}(A)}_{*-i+1,\bu+i}\}$$ are related by
$$\sum_{i+j=k}d_1^{\,i}s^{j-1}+s^{j-1}d_1^{\,i}=\left\{\begin{array}{ll}
1_{\,{^\varrho{\cal M}(A)}},&k=1,\\ 0,&k\not=1.
\end{array}\right.\eqno(3.6)$$

Given a stable $\D$-module $(^\varrho{\cal M}(A),d_1^{\,i})$,
consider the corresponding chain complex $(\overline{^\varrho{\cal
M}(A)},b^{'})$, where
$b^{'}=(d_1^{\,0}+d_1^{\,1}+\dots+d_1^{\,i}+\dots):\overline{^\varrho{\cal
M}(A)}_\bu\to \overline{^\varrho{\cal M}(A)}_{\bu-1}$. It follows
from $(3.6)$ that the map
$$s=(s^{-1}+s^0+s^1+\dots+s^i+\dots\,):\overline{^\varrho{\cal
M}(A)}_\bu\to \overline{^\varrho{\cal M}(A)}_{\bu+1}$$ satisfies
the condition $b^{'}s+sb^{'}=1_{\overline{^\varrho{\cal
M}(A)}}\,$, which means that $s$ is a contracting homotopy for the
chain complex $(\overline{^\varrho{\cal M}(A)},b^{'})$.

Given an arbitrary involutive $\A^{hu}$-algebra
$(A,d,\pi_n,\tau_n^{j_q,\dots,j_1},*)$, the relations $(1.9)$
follows that there is a left action of the group
$\mathbb{Z}_2=\{1,\vartheta\,|\,\vartheta^2=1\}$ on the chain
complex $(\overline{^\varrho{\cal M}(A)},b^{'})$. This left action
is defined by means of the module automorphism
$\vartheta:\overline{^\varrho{\cal
M}(A)}\to\overline{^\varrho{\cal M}(A)}$ of the order two, which
at any element $x\in \overline{^\varrho{\cal M}(A)}_n$ is given by
$\vartheta(x)=\overline{R}_n\overline{T}_n(x)$.

Consider the chain bicomplex $(Q(\overline{^\varrho{\cal
M}(A)}),-b^{'},D')$ defined by $$Q(\overline{^\varrho{\cal
M}(A)})_{n,m}=\overline{^\varrho{\cal M}(A)}_n,\quad n\geqslant
0,\quad m\geqslant 0,$$
$$D'=(-1)^{n+1}(1+(-1)^{m+1}\overline{R}_n\overline{T}_n):D':Q(\overline{^\varrho{\cal
M}(A)})_{n,m}\to Q(\overline{^\varrho{\cal M}(A)})_{n,m-1},$$
$$-b^{'}:Q(\overline{^\varrho{\cal
M}(A)})_{n,m}=\overline{^\varrho{\cal M}(A)}_n\to
\overline{^\varrho{\cal M}(A)}_{n-1}=Q(\overline{^\varrho{\cal
M}(A)})_{n-1,m}.$$ For the chain complex associated with the chain
complex $(Q(\overline{^\varrho{\cal M}(A)}),-b^{'},D')$, we use
the notation $({\rm Tot}(Q(\overline{^\varrho{\cal
M}(A)})),\widehat{D}')$, where $\widehat{D}'=-b+D'$. If we
consider the spectral sequence of the chain bicomplex
$(Q(\overline{^\varrho{\cal M}(A)}),-b^{'},D')$ and use the the
contractibility of the chain complex $(\overline{^\varrho{\cal
M}(A)},b^{'})$, then we obtain $H_n({\rm
Tot}(Q(\overline{^\varrho{\cal M}(A)})),\widehat{D}')=0$ for all
$n\geqslant 0$.

In what follows, by the dihedral homology $^\varrho HD(A)$ of any
involutive $\A^{hu}$-algebra
$(A,d,\pi_n,\tau_n^{j_q,\dots,j_1},*)$ we mean the dihedral
homology $^\varrho HD(A)$ of the corresponding involutive
$\A$-algebra $(A,d,\pi_n,*)$. Similarly, by the reflexive homology
$^\varrho HR(A)$ of any involutive $\A^{hu}$-algebra
$(A,d,\pi_n,\tau_n^{j_q,\dots,j_1},*)$ we mean the reflexive
homology $^\varrho HR(A)$ of the corresponding $\A$-algebra
$(A,d,\pi_n,*)$.

Now, using the acyclicity of the chain complex $({\rm
Tot}(Q(\overline{^\varrho{\cal M}(A)})),\widehat{D}')$, we find a
relationship between the dihedral and the reflexive homology of
involutive homotopy unital $\A$-algebras.

{\bf Theorem 3.1}. For any involutive $\A^{hu}$-algebra
$(A,d,\pi_n,\tau_n^{j_q,\dots,j_1},*)$, there is the long exact
sequence $$\cdots\buildrel{\delta_*}\over\longrightarrow\,\!
^\varrho HR_n(A)\buildrel{i_*}\over\longrightarrow\,\!^\varrho
HD_n(A)\buildrel{p_*}\over\longrightarrow\,\!\!
^{-\varrho}HD_{n-2}(A)\buildrel{\delta_*}\over\longrightarrow\,\!^\varrho
HR_{n-1}(A)\buildrel{i_*}\over\longrightarrow\cdots\,,\eqno(3.7)$$
where maps $i_*$ and $p_*$ are induced respectively by the
inclusion and by the projection, and the map $\delta_*$ is induced
by the connecting homomorphism.

{\bf Proof}. Given the triple chain complex
$(D(\overline{^\varrho{\cal M}(A)}),\delta_1,\delta_2,\delta_3)$,
consider the chain bicomplex
$(\widetilde{D}(\overline{^\varrho{\cal
M}(A)}),\widetilde{\delta}_1,\widetilde{\delta}_2)$ defined by
$$\widetilde{D}(\overline{^\varrho{\cal
M}(A)})_{n,m}=\bigoplus_{k+l=n}D(\overline{^\varrho{\cal
M}(A)})_{k,m,l},\quad n\geqslant 0,\quad m\geqslant 0,\quad
\widetilde{\delta}_1=\delta_1+\delta_3,\quad
\widetilde{\delta}_2=\delta_2.$$ The chain complex associated with
the chain bicomplex $(\widetilde{D}(\overline{^\varrho{\cal
M}(A)}),\widetilde{\delta}_1,\widetilde{\delta}_2)$ we denote by
$({\rm Tot}(\widetilde{D}(\overline{^\varrho{\cal
M}(A)})),\widetilde{\delta})$, where
$\widetilde{\delta}=\widetilde{\delta}_1+\widetilde{\delta}_2$. It
is easy to see that the chain complex $({\rm
Tot}(D(\overline{^\varrho{\cal M}(A)})),\delta)$ is isomorphic to
the chain complex $({\rm
Tot}(\widetilde{D}(\overline{^\varrho{\cal
M}(A)})),\widetilde{\delta})$ and, consequently, we have the
isomorphism $H({\rm Tot}(\widetilde{D}(\overline{^\varrho{\cal
M}(A)})),\widetilde{\delta})=\,\!^\varrho HD(A)$. Consider in the
chain complex $({\rm Tot}(\widetilde{D}(\overline{^\varrho{\cal
M}(A)})),\widetilde{\delta})$ the chain subcomplex
$(P(\overline{^\varrho{\cal M}(A)}),\widetilde{\delta})$ defined
by the following formulae: $$P(\overline{^\varrho{\cal
M}(A)})_n=\widetilde{D}(\overline{^\varrho{\cal
M}(A)})_{n-1,1}\oplus \widetilde{D}(\overline{^\varrho{\cal
M}(A)})_{n,0},\quad n\geqslant 0,$$
$$\widetilde{\delta}(x_{n-1},x_n)=(\widetilde{\delta}_1(x_{n-1}),\widetilde{\delta}_1(x_n)
+\widetilde{\delta}_2(x_{n-1}),\quad (x_{n-1},x_n)\in
P(\overline{^\varrho{\cal M}(A)})_n.$$ It is easy to see that
there are the equalities of modules
$$\widetilde{D}(\overline{^\varrho{\cal M}(A)})_{n-1,1}={\rm
Tot}(Q(\overline{^\varrho{\cal M}(A)}))_{n-1},\quad
\widetilde{D}(\overline{^\varrho{\cal M}(A)})_{n,0}={\rm
Tot}(R(\overline{^\varrho{\cal M}(A)}))_n.$$ Consider the short
exact sequence of chain complexes $$0\to
(P(\overline{^\varrho{\cal
M}(A)}),\widetilde{\delta})\buildrel{j}\over\longrightarrow ({\rm
Tot}(\widetilde{D}(\overline{^\varrho{\cal
M}(A)})),\widetilde{\delta})\buildrel{p}\over\longrightarrow
\Sigma^{-2}({\rm Tot}(\widetilde{D}(\overline{^\varrho{\cal
M}(A)})),\widetilde{\delta})\to 0,$$
$$j((x_{n-1},x_n))=(0,\dots,0,x_{n-1},x_n),\quad n\geqslant 0,$$
$$\Sigma^{-2}({\rm Tot}(\widetilde{D}(\overline{^\varrho{\cal
M}(A)}))_\bu,\widetilde{\delta})=({\rm
Tot}(\widetilde{D}(\overline{^{-\varrho}{\cal
M}(A)}))_{\bu-2},\widetilde{\delta}),$$
$$p((x_0,\dots,x_{n-2},x_{n-1},x_n))=(x_0,\dots,x_{n-2}),\quad
n\geqslant 0.$$ In the homology, this short exact sequence induces
the long exact sequence
$$\cdots\buildrel{\delta}\over\longrightarrow
H_n(P(\overline{^\varrho{\cal
M}(A)}))\buildrel{j_*}\over\longrightarrow\,\!^\varrho
HD_n(A)\buildrel{p_*}\over\longrightarrow\,\!\!$$
$$^{-\varrho}HD_{n-2}(A)\buildrel{\delta}\over\longrightarrow
H_{n-1}(P(\overline{^\varrho{\cal
M}(A)}))\buildrel{i_*}\over\longrightarrow\cdots\,.\eqno(3.8)$$
Now, let us prove that the graded modules
$H(P(\overline{^\varrho{\cal M}(A)}))$ and $^\varrho HR(A)$ are
isomorphic. To this purpose we consider the short exact sequence
of chain complexes $$0\to ({\rm Tot}(R(\overline{^\varrho{\cal
M}(A)}))_\bu,\widehat{D})\buildrel{\alpha}\over\longrightarrow
(P(\overline{^\varrho{\cal
M}(A)})_\bu,\widetilde{\delta})\buildrel{\beta}\over\longrightarrow
({\rm Tot}(Q(\overline{^\varrho{\cal
M}(A)}))_{\bu-1},\widehat{D}')\to 0,$$ where
$\alpha(x_n)=(0,x_n)$, $\beta(x_{n-1},x_n)=x_{n-1}$,
$(x_{n-1},x_n)\in (P(\overline{^\varrho{\cal M}(A)})_n$. In the
homology, this short exact sequence induces the long exact
sequence. Since the equality $H_{\bu-1}({\rm
Tot}(Q(\overline{^\varrho{\cal M}(A)})),\widehat{D}')=0$ is true,
this long exact sequence implies that the map
$\alpha_*:\,\!^\varrho HR_\bu(A)\to
H_\bu(P(\overline{^\varrho{\cal M}(A)}))$ is an isomorphism. Now,
replacing the modules $H_\bu(P(\overline{^\varrho{\cal M}(A)}))$
in the long exact sequence $(3.8)$ by $^\varrho HR_\bu(A)$ and the
maps $j_*$ and $\delta$ by $i_*=j_*\alpha_*$ and
$\delta_*=\alpha^{-1}_*\delta$, respectively, we obtain the long
exact sequence $(3.7)$.~~~$\blacksquare$

In conclusion we mention that if an involutive $\A^{hu}$-algebra
$(A,d,\pi_n,\tau_n^{j_q,\dots,j_1},*)$ is a differential
involutive unital associative algebra $(A,d,\pi,u,*)$, where
$\pi=\pi_0$, $u=\tau_0^0$ and $\tau_n^{j_q,\dots,j_1}=0$ in all
other cases, then the exact sequence $(3.7)$ coincides with the
well-known Krasauskas-Lapin-Solov'ev exact sequence that
connecting the dihedral and the reflexive homology of an algebra
$(A,d,\pi,u,*)$ over any unital commutative ring.

\vspace{1cm}

Saransk, Russia,

e-mail: slapin@mail.ru

\end{document}